\documentclass[11pt]{amsart}
\usepackage{amsmath,amssymb,amsthm,stmaryrd}
\newtheorem{theorem}{Theorem}[section]
\newtheorem{proposition}[theorem]{Proposition}
\newtheorem{lemma}[theorem]{Lemma}
\newtheorem{corollary}[theorem]{Corollary}
\theoremstyle{definition}
\newtheorem{definition}[theorem]{Definition}
\newtheorem{remark}[theorem]{Remark}

\begin{document}

\title{Ricci flow with surgery in higher dimensions}
\author{Simon Brendle}
\address{Department of Mathematics \\ Columbia University \\ New York, NY 10027}
\thanks{The author is grateful to Professor Hong Huang for comments on an earlier version of this paper. This project was supported in part by the National Science Foundation under grants DMS-1201924 and DMS-1505724. The author would like to acknowledge the hospitality of T\"ubingen University, where part of this work was carried out.}
\begin{abstract}
We present a new curvature condition which is preserved by the Ricci flow in higher dimensions. For initial metrics satisfying this condition, we establish a higher dimensional version of Hamilton's neck-like curvature pinching estimate. Using this estimate, we are able to prove a version of Perelman's Canonical Neighborhood Theorem in higher dimensions. This makes it possible to extend the flow beyond singularities by a surgery procedure in the spirit of Hamilton and Perelman. As a corollary, we obtain a classification of all diffeomorphism types of such manifolds in terms of a connected sum decomposition. In particular, the underlying manifold cannot be an exotic sphere.

Our result is sharp in many interesting situations. For example, the curvature tensors of $\mathbb{CP}^{n/2}$, $\mathbb{HP}^{n/4}$, $S^{n-k} \times S^k$ ($2 \leq k \leq n-2$), $S^{n-2} \times \mathbb{H}^2$, $S^{n-2} \times \mathbb{R}^2$ all lie on the boundary of our curvature cone. Another borderline case is the pseudo-cylinder: this is a rotationally symmetric hypersurface which is weakly, but not strictly, two-convex. Finally, the curvature tensor of $S^{n-1} \times \mathbb{R}$ lies in the interior of our curvature cone.
\end{abstract}
\maketitle

\section{Introduction}

Since its introduction by Hamilton \cite{Hamilton1} in 1982, the Ricci flow has become a fundamental tool in Riemannian geometry. In particular, two lines of research have been pursued: 

First, there are a number of theorems which show, under suitable assumptions on the initial metric, that the Ricci flow converges to a metric of constant curvature, up to rescaling. The earliest result in this direction is the famous work of convergence theorem of Hamilton \cite{Hamilton1} for three-manifolds with positive Ricci curvature. In higher dimensions, Huisken \cite{Huisken}, Margerin \cite{Margerin1}, and Nishikawa \cite{Nishikawa} found various pinching conditions that guarantee that the flow converges to a round metric, after rescaling. These conditions were later weakened in work of Margerin \cite{Margerin2}, \cite{Margerin3}. In a fundamental work \cite{Hamilton2}, Hamilton introduced his PDE-ODE principle, and used it to show that the Ricci flow deforms any four-manifold with positive curvature operator to a round metric. This was generalized to higher dimensions in an important paper by B\"ohm and Wilking \cite{Bohm-Wilking}. In \cite{Brendle-Schoen}, the author and R.~Schoen showed that, if the initial metric has $1/4$-pinched sectional curvature, then the Ricci flow will converge to a round metric after rescaling. In particular, this result implies the Differentiable Sphere Theorem: a Riemannian manifold with $1/4$-pinched sectional curvature is diffeomorphic to a space form. This result was further generalized in \cite{Brendle1}.

Second, it is of interest to find conditions that restrict the singularities that can form under the evolution to so-called neck pinch singularities. This was first done in a seminal paper by Hamilton \cite{Hamilton5}, where he used the Ricci flow to classify four-manifolds with positive isotropic curvature. In a striking breakthrough, Perelman \cite{Perelman1},\cite{Perelman2},\cite{Perelman3} succeeded in carrying out a similar program in dimension $3$, without any assumptions on the initial metric, proving the Poincar\'e conjecture as a direct consequence. 

In order to understand the global behavior of the Ricci, it is important to find curvature conditions which are preserved by the Ricci flow. For example, Hamilton proved that the positive scalar curvature and positive curvature operator are preserved by the flow. In \cite{Brendle-Schoen}, the author and R.~Schoen found additional curvature conditions that are preserved by the Ricci flow. These include the condition of positive isotropic curvature (PIC), which originated in the work of Micallef and Moore \cite{Micallef-Moore} in minimal surface theory, as well as the PIC1 and PIC2 conditions. By definition, $M$ satisfies the PIC1 condition if $M \times \mathbb{R}$ has nonnegative isotropic curvature, and $M$ satisfies the PIC2 condition if $M \times \mathbb{R}^2$ has nonnegative isotropic curvature. Note that the PIC condition is weaker than the PIC1 condition, the PIC1 condition is weaker than the PIC2 condition, and the PIC2 condition is weaker than positivity of the curvature operator or $1/4$-pinching. 

Our goal in this paper is to prove a higher-dimensional version of Hamilton's theorem on Ricci flow with surgery \cite{Hamilton5}. To that end, we construct new curvature conditions that are preserved by the Ricci flow, and which pinch toward a better curvature condition. To find such conditions, we need to find convex subsets of the space of algebraic curvature tensors which are preserved by the ODE $\frac{d}{dt} R = Q(R)$, where $Q(R)$ is a suitable quadratic polynomial in the curvature tensor. We will decompose an algebraic curvature tensor as a sum $R = S + H \owedge \text{\rm id}$, where $H$ is a symmetric bilinear form, $S$ is a curvature tensor with $\text{\rm Ric}_0(S)=0$, and $\owedge$ denotes the Kulkarni-Nomizu product. In other words, $H$ determines the trace-free tensor of $R$, while $S$ determines the Weyl tensor of $R$. Note that the decomposition is not unique, as we are free to add a multiple of the identity to $H$ and subtract a multiple of the curvature tensor of $S^n$ from $S$. We can now evolve $S$ and $H$ by a coupled system of ODEs, and look for conditions on $S$ and $H$ that are preserved by this ODE system; these conditions combine a positivity condition on $S$ with an upper bound for the largest eigenvalue of $H$. This leads to new preserved curvature conditions which allow a surgery construction in higher dimensions:

\begin{theorem}
\label{surgically.modified.flow}
Let $(M,g_0)$ be a compact manifold of dimension $n \geq 5$ whose curvature tensor lies in the interior of the cone 
\begin{align*} 
&\{R = S + H \owedge \text{\rm id}: S \in PIC2, \, \text{\rm Ric}_0(S)=0, \\ 
&\hspace{32mm} \text{\rm tr}(H) \, \text{\rm id} - (n-4) \, H \geq 0\} 
\end{align*} 
at each point. Then the curvature tensor remains in this set if the metric is evolved by the Ricci flow. Moreover, if $M$ does not contain non-trivial incompressible space forms $S^{n-1}/\Gamma$, then there exists a Ricci flow with surgery starting from $(M,g_0)$. This flow involves performing finitely many surgeries on necks of the form $S^{n-1} \times I$. Finally, the surgically modified flow becomes extinct in finite time. 
\end{theorem}

\begin{theorem}
\label{classification}
Under the assumptions of Theorem \ref{surgically.modified.flow}, $M$ is diffeomorphic to a connected sum of finitely many pieces. Each piece is a quotient of $S^n$ or a compact quotient of $S^{n-1} \times \mathbb{R}$ by standard isometries. 
\end{theorem}

\begin{remark}
\begin{itemize}
\item Up to scaling, the curvature tensor of $S^{n-1} \times \mathbb{R}$ can be expressed as 
\[\begin{bmatrix} 0 & & & \\ & 1 & & \\ & & \ddots & \\ & & & 1 \end{bmatrix} \owedge \begin{bmatrix} 0 & & & \\ & 1 & & \\ & & \ddots & \\ & & & 1 \end{bmatrix} = \begin{bmatrix} -1 & & & \\ & 1 & & \\ & & \ddots & \\ & & & 1 \end{bmatrix} \owedge \text{\rm id}.\]
Therefore, the curvature tensor of $S^{n-1} \times \mathbb{R}$ lies in the interior of our curvature cone. Consequently, our curvature conditions is preserved under formation of connected sums. 
\item Consider next a pseudo-cylinder, that is, a rotationally symmetric hypersurface which is weakly, but not strictly, two-convex. Up to scaling, the curvature tensor of a pseudo-cylinder is given by 
\[\begin{bmatrix} -1 & & & \\ & 1 & & \\ & & \ddots & \\ & & & 1 \end{bmatrix} \owedge \begin{bmatrix} -1 & & & \\ & 1 & & \\ & & \ddots & \\ & & & 1 \end{bmatrix} = \begin{bmatrix} -3 & & & \\ & 1 & & \\ & & \ddots & \\ & & & 1 \end{bmatrix} \owedge \text{\rm id}.\] 
This curvature tensor lies on the boundary of our curvature cone. Note that $S=0$ in this example.
\item The curvature tensor of $S^{n-2} \times \mathbb{H}^2$ lies on the boundary of our curvature cone. In this example, $S=0$.
\item The set 
\[\{H \owedge \text{\rm id}: \text{\rm tr}(H) \, \text{\rm id} - (n-4) \, H \geq 0\}\] 
can be interpreted as the convex hull of the set of all curvature tensors of pseudo-cylinders. 
\item The curvature tensors of $\mathbb{CP}^{n/2}$, $\mathbb{HP}^{n/4}$, and $S^{n-k} \times S^k$ ($2 \leq k \leq n-2$) all lie on the boundary of our curvature cone. In each of these examples, $H=0$.
\item The curvature tensor of $S^{n-2} \times \mathbb{R}^2$ lies on the boundary of our curvature cone, since it can be expressed as a sum of the curvature tensor of $S^{n-2} \times S^2$ and the curvature tensor of $S^{n-2} \times \mathbb{H}^2$.
\item Theorem \ref{surgically.modified.flow} also holds for in dimension $4$. For $n=4$, any curvature tensor satisfying $a_1 > 0$ and $c_1 > 0$ (in the notation of \cite{Hamilton2} and \cite{Hamilton5}) lies in our curvature cone.
\end{itemize}
\end{remark}

\begin{remark}
As in \cite{Hamilton5}, the assumption that $M$ does not contain non-trivial incompressible $(n-1)$-dimensional space forms is used to rule out quotient necks which are modeled on noncompact quotients $(S^{n-1} \times \mathbb{R})/\Gamma$; a precise statement can be found in the appendix. If we remove this assumption, quotient necks may occur. Doing surgery on a quotient neck will produce an orbifold. Thus, it becomes necessary to work in the class of orbifolds. For $n=4$, the necessary adaptations in the orbifold setting are discussed by Chen, Tang, and Zhu \cite{Chen-Tang-Zhu}. We expect that the arguments in this paper can similarly be extended to the orbifold setting; the details will appear elsewhere.
\end{remark}

In Section \ref{pinching}, we construct new preserved curvature conditions in higher dimensions. In particular, we show that, under the Hamilton ODE $\frac{d}{dt} R = Q(R)$, the curvature cone in Theorem \ref{surgically.modified.flow} pinches toward a smaller cone which is contained in the PIC2 cone. This can be viewed as a cylindrical estimate for the Ricci flow in higher dimensions. 

In Section \ref{splitting.theorems}, we prove various splitting theorems, which will be needed in the later sections. 

In Section \ref{ancient.solutions}, we study ancient $\kappa$-solutions which satisfy the cylindrical estimate in Section \ref{pinching}. We prove an analogue of Perelman's long range curvature estimate and Perelman's universal noncollapsing theorem in this setting. As a consequence, we obtain an analogue of Perelman's Canonical Neighborhood Theorem for ancient $\kappa$-solutions.

In Section \ref{canonical.neighborhood.theorem}, we study a solution to the Ricci flow on a compact manifold, where the initial metric satisfies the curvature condition in Theorem \ref{surgically.modified.flow}. We prove that the high curvature regions are modeled on ancient $\kappa$-solutions. In particular, this gives a Canonical Neighborhood Theorem which holds at all points in space-time where the curvature is sufficiently large.

Finally, in Section \ref{first.singular.time}, we give a precise description of the behavior of the flow at the first singular time. Moreover, we verify that our curvature pinching estimates are preserved under surgery. 

\section{Curvature pinching estimates for the Ricci flow in higher dimensions}

\label{pinching}

Given two algebraic curvature tensors $S$ and $T$, we define a new algebraic curvature tensor $B(S,T)$ by 
\begin{align*} 
B(S,T)_{ijkl} 
&= \frac{1}{2} \sum_{p,q=1}^n (S_{ijpq} \, T_{klpq} + S_{klpq} \, T_{ijpq}) \\ 
&+ \sum_{p,q=1}^n (S_{ipkq} \, T_{jplq} - S_{iplq} \, T_{jpkq} - S_{jpkq} \, T_{iplq} + S_{jplq} \, T_{ipkq})
\end{align*} 
Clearly, $B(S,T) = B(T,S)$. Moreover, $B(R,R) = Q(R)$, where $Q(R)$ is the term appearing in Hamilton's curvature ODE. Given an algebraic curvature tensor $S$ and a symmetric bilinear form $H$, we define 
\[(S*H)_{ik} = \sum_{j,l=1}^n S_{ijkl} \, H_{jl}.\] 
Finally, we denote by $(A \owedge B)_{ijkl} = A_{ik} B_{jl} - A_{il} B_{jk} - A_{jk} B_{il} + A_{jl} B_{ik}$ the Kulkarni-Nomizu product of two symmetric two-tensors $A$ and $B$.

\begin{lemma}
\label{algebra}
We have 
\[B(S,H \owedge \text{\rm id}) = \text{\rm Ric}(S) \owedge H  + (S*H) \owedge \text{\rm id}.\] 
and 
\begin{align*} 
Q(H \owedge \text{\rm id}) 
&= (n-2) \, H \owedge H + 2 \, \text{\rm tr}(H) \, H \owedge \text{\rm id} \\ 
&- 2 \, H^2 \owedge \text{\rm id} + |H|^2 \, \text{\rm id} \owedge \text{\rm id}. 
\end{align*} 
\end{lemma}

\textbf{Proof.} 
We begin with the first statement. If we put $T = H \owedge \text{\rm id}$, then we obtain 
\begin{align*} 
\frac{1}{2} \sum_{p,q=1}^n S_{ijpq} \, T_{klpq} 
&= \frac{1}{2} \sum_{p,q=1}^n S_{ijpq} \, (H_{kp} \delta_{lq} - H_{kq} \delta_{lp} - H_{lp} \delta_{kq} + H_{lq} \delta_{kp}) \\ 
&= \sum_{p=1}^n S_{ijpl} \, H_{kp} + \sum_{p=1}^n S_{ijkp} \, H_{lp} 
\end{align*} 
and 
\begin{align*} 
\sum_{p,q=1}^n S_{ipkq} \, T_{jplq} 
&= \sum_{p,q=1}^n S_{ipkq} \, (H_{jl} \delta_{pq} - H_{jq} \delta_{lp} - H_{lp} \delta_{jq} + H_{pq} \delta_{jl}) \\ 
&= \text{\rm Ric}(S)_{ik} \, H_{jl} - \sum_{p=1}^n S_{ilkp} \, H_{jp} - \sum_{p=1}^n S_{ipkj} \, H_{lp} + (S*H)_{ik} \, \delta_{jl}. 
\end{align*} 
Putting these facts together, we obtain 
\begin{align*} 
&B(S,T)_{ijkl} - (\text{\rm Ric}(S) \owedge H)_{ijkl} - ((S*H) \owedge \text{\rm id})_{ijkl} \\ 
&= \sum_{p=1}^n S_{ijpl} \, H_{kp} + \sum_{p=1}^n S_{ijkp} \, H_{lp} \\ 
&+ \sum_{p=1}^n S_{pjkl} \, H_{ip} + \sum_{p=1}^n S_{ipkl} \, H_{jp} \\ 
&- \sum_{p=1}^n S_{ilkp} \, H_{jp} - \sum_{p=1}^n S_{ipkj} \, H_{lp} \\ 
&+ \sum_{p=1}^n S_{iklp} \, H_{jp} + \sum_{p=1}^n S_{iplj} \, H_{kp} \\ 
&+ \sum_{p=1}^n S_{jlkp} \, H_{ip} + \sum_{p=1}^n S_{jpki} \, H_{lp} \\ 
&- \sum_{p=1}^n S_{jklp} \, H_{ip} - \sum_{p=1}^n S_{jpli} \, H_{kp}, 
\end{align*}
and the right hand side vanishes by the first Bianchi identity. This proves the first statement.

To derive the second statement, we apply the first statement with $S = H \owedge \text{\rm id}$. Then $\text{\rm Ric}(S) = (n-2) \, H + \text{\rm tr}(H) \, \text{\rm id}$. Moreover, $S*H = \text{\rm tr}(H) \, H - 2 \, H^2 + |H|^2 \, \text{\rm id}$. From this, the assertion follows. \\

\begin{definition}
Given $\sigma \in (0,2]$ and $\theta \geq 0$, we define a cone $\mathcal{C}_{\sigma,\theta}$ in the space of algebraic curvature tensors by 
\begin{align*} 
\mathcal{C}_{\sigma,\theta} := &\{R = S + H \owedge \text{\rm id}: S \in PIC2, \, \text{\rm Ric}_0(S)=0, \\ 
&\hspace{32mm} \text{\rm tr}(H) \, \text{\rm id} - (n-2\sigma) \, H \geq 0, \\ 
&\hspace{32mm} \text{\rm tr}(H) - \theta \, \text{\rm scal}(S) \geq 0\}. 
\end{align*} 
Equivalently, $\mathcal{C}_{\sigma,\theta}$ consists of all algebraic curvature tensors $R$ with the property that 
\[R - \frac{1}{n-2} \, \text{\rm Ric}_0 \owedge \text{\rm id} - \frac{1}{n} \, \frac{\theta}{1+2(n-1)\theta} \, \text{\rm scal} \, \text{\rm id} \owedge \text{\rm id} \in PIC2\] 
and 
\[|v|^2 \, R - \frac{1}{n-2} \, |v|^2 \, \text{\rm Ric}_0 \owedge \text{\rm id} - \frac{n-2\sigma}{2(n-2)\sigma} \, \text{\rm Ric}_0(v,v) \, \text{\rm id} \owedge \text{\rm id} \in PIC2\] 
for each $v$.
\end{definition}

We first examine the properties of the cones $\mathcal{C}_{1,0}$ and $\mathcal{C}_{2,0}$.

\begin{proposition} 
\label{C1}
Suppose that $R \in \mathcal{C}_{1,0}$. Then $R \in PIC2$. Moreover, if $\text{\rm Ric}(v,v) = 0$ for some unit vector $v$, then $R = c \, (\text{\rm id} - 2 \, v \otimes v) \owedge \text{\rm id}$ for some $c \geq 0$. In other words, either $R=0$ or $R$ is the curvature tensor of a cylinder $S^{n-1} \times \mathbb{R}$. 
\end{proposition}

\textbf{Proof.} 
By definition, we may write $R = S + H \owedge \text{\rm id}$, where $S \in PIC2$, $\text{\rm Ric}_0(S)=0$, and $\text{\rm tr}(H) \, \text{\rm id} - (n-2) \, H \geq 0$. The condition $\text{\rm tr}(H) \, \text{\rm id} - (n-2) \, H \geq 0$ easily implies that $H$ is weakly two-positive. Consequently, $H \owedge \text{\rm id} \in PIC2$, hence $R \in PIC2$. Moreover, the Ricci tensor of $R$ is given by $\frac{1}{n} \, \text{\rm scal}(S) \, \text{\rm id} + \text{\rm tr}(H) \, \text{\rm id} + (n-2) \, H$. The condition $\text{\rm tr}(H) \, \text{\rm id} - (n-2) \, H \geq 0$ implies $\text{\rm tr}(H) \, \text{\rm id} + (n-2) \, H \geq 0$. Hence, if $\text{\rm Ric}(v,v) = 0$, then we have $\text{\rm scal}(S) = 0$ and furthermore $H(v,v) = -\frac{1}{n-2} \, \text{\rm tr}(H)$. From this, we deduce that $S = 0$ and $H = c \, (\text{\rm id} - 2 \, v \otimes v)$ for some $c \geq 0$. \\

\begin{proposition} 
\label{C2}
The curvature tensor of $S^{n-1} \times \mathbb{R}$ lies in the interior of the cone $\mathcal{C}_{2,0}$. Moreover, the curvature tensors of $S^{n-2} \times \mathbb{R}^2$ and $S^{n-2} \times S^2$ lie in $\mathcal{C}_{2,0}$. Here, the curvature tensor of $S^{n-2} \times S^2$ is normalized so that the trace-free Ricci part vanishes. Finally, the curvature tensor of a pseudo-cylinder lies in $\mathcal{C}_{2,0}$. 
\end{proposition}

\textbf{Proof.} Let 
\[S_{ijkl} = \begin{cases} (n-3) \, (\delta_{ik} \delta_{jl} - \delta_{il} \delta_{jk}) & \text{\rm if $i,j,k,l \in \{1,2\}$} \\ \delta_{ik} \delta_{jl} - \delta_{il} \delta_{jk} & \text{\rm if $i,j,k,l \in \{3,\hdots,n\}$} \\ 0 & \text{\rm otherwise.} \end{cases}\] 
Then $S \in PIC2$, and the trace-free Ricci part of $S$ vanishes. Geometrically, $S$ is the curvature tensor of $S^2 \times S^{n-2}$ (suitably normalized). Moreover, let 
\[H_{ij} = \begin{cases} -\delta_{ij} & \text{\rm if $i,j \in \{1,2\}$} \\ \delta_{ij} & \text{\rm if $i,j \in \{3,\hdots,n\}$} \\ 0 & \text{\rm otherwise.} \end{cases}\] 
Clearly, $\text{\rm tr}(H) \, \text{\rm id} - (n-4) \, H \geq 0$. Geometrically, $H \owedge \text{\rm id}$ represents the curvature tensor of $\mathbb{H}^2 \times S^{n-2}$ (suitably normalized). Finally, we observe that $S + \frac{n-3}{2} \, H \owedge \text{\rm id}$ is the curvature tensor of $\mathbb{R}^2 \times S^{n-2}$. This shows that the curvature tensors of $S^2 \times S^{n-2}$ and $\mathbb{R}^2 \times S^{n-2}$ both lie in the cone $\mathcal{C}_{2,0}$. Finally, the fact that the pseudo-cylinder lies in $\mathcal{C}_{2,0}$ is clear from the definition. \\

We now state the main result of this section:

\begin{theorem}
\label{pinching.cones}
For each $n$, there exists a positive constant $\bar{\theta} = \bar{\theta}(n)$ with the following property: For each $\sigma \in (0,2]$ and each $\theta \in [0,\bar{\theta}]$, the cone $\mathcal{C}_{\sigma,\theta}$ is invariant under the Hamilton ODE. Moreover, if $\sigma \in (0,1) \cup (1,2)$ and $\theta \in (0,\bar{\theta})$, the cone $\mathcal{C}_{\sigma,\theta}$ is transversally invariant away from $0$.
\end{theorem}

\textbf{Proof.} 
In the following, we fix real numbers $\sigma \in (0,2]$ and $\theta \geq 0$. We evolve $S$ and $H$ by 
\begin{align*} 
\frac{d}{dt} S 
&= Q(S) + (n-2) \, H \owedge H - 2 \, \text{\rm tr}(H) \, H \owedge \text{\rm id} + 2 \, H^2 \owedge \text{\rm id} \\ 
&+ \frac{2}{\sigma(n-2\sigma)} \, \text{\rm tr}(H)^2 \, \text{\rm id} \owedge \text{\rm id} - \frac{2-\sigma}{\sigma} \, |H|^2 \, \text{\rm id} \owedge \text{\rm id}
\end{align*}
and 
\begin{align*} 
\frac{d}{dt} H 
&= 2 \, S*H + \frac{2}{n} \, \text{\rm scal}(S) \, H + 4 \, \text{\rm tr}(H) \, H - 4 \, H^2 \\ 
&- \frac{2}{\sigma(n-2\sigma)} \, \text{\rm tr}(H)^2 \, \text{\rm id} + \frac{2}{\sigma} \, |H|^2 \, \text{\rm id}. 
\end{align*}
We proceed in several steps: \\

\textit{Step 1:} We claim that the condition $\text{\rm Ric}_0(S) = 0$ is preserved. To show this, it suffices to prove that the term 
\begin{align*} 
T &:= (n-2) \, H \owedge H - 2 \, \text{\rm tr}(H) \, H \owedge \text{\rm id} + 2 \, H^2 \owedge \text{\rm id} \\ 
&+ \frac{2}{\sigma(n-2\sigma)} \, \text{\rm tr}(H)^2 \, \text{\rm id} \owedge \text{\rm id} - \frac{2-\sigma}{\sigma} \, |H|^2 \, \text{\rm id} \owedge \text{\rm id} 
\end{align*} 
has vanishing trace-free Ricci part. 

The Ricci tensor of $(n-2) \, H \owedge H - 2 \, \text{\rm tr}(H) \, H \owedge \text{\rm id} + 2 \, H^2 \owedge \text{\rm id}$ is given by 
\begin{align*} 
&(n-2) \, [2 \, \text{\rm tr}(H) \, H - 2 \, H^2] \\ 
&- 2 \, [(n-2) \, \text{\rm tr}(H) \, H + \text{\rm tr}(H)^2 \, \text{\rm id}] \\ 
&+ 2 \, [(n-2) \, H^2 + |H|^2 \, \text{\rm id}] \\ 
&= 2 \, [|H|^2 - \text{\rm tr}(H)^2] \, \text{\rm id}. 
\end{align*}
This shows that the trace-free Ricci part of $T$ vanishes. \\

\textit{Step 2:} We claim that the condition $S \in PIC2$ is preserved. It suffices to show that the term 
\begin{align*} 
T &:= (n-2) \, H \owedge H - 2 \, \text{\rm tr}(H) \, H \owedge \text{\rm id} + 2 \, H^2 \owedge \text{\rm id} \\ 
&+ \frac{2}{\sigma(n-2\sigma)} \, \text{\rm tr}(H)^2 \, \text{\rm id} \owedge \text{\rm id} - \frac{2-\sigma}{\sigma} \, |H|^2 \, \text{\rm id} \owedge \text{\rm id} 
\end{align*} 
has nonnegative curvature operator. A straightforward calculation gives 
\begin{align*} 
T &= (n-2) \, A \owedge A + 2 \, A^2 \owedge \text{\rm id} - 2 \, \text{\rm tr}(A) \, A \owedge \text{\rm id} \\ 
&+ \frac{1}{\sigma^2} \, \text{\rm tr}(A)^2 \, \text{\rm id} \owedge \text{\rm id} - \frac{2-\sigma}{\sigma} \, |A|^2 \, \text{\rm id} \owedge \text{\rm id}, 
\end{align*} 
where $A := \frac{1}{n-2\sigma} \, \text{\rm tr}(H) \, \text{\rm id} -  H \geq 0$. To show that $T$ has nonnegative curvature operator, it therefore suffices to show that 
\begin{align*} 
&(n-2) \, a_i a_j + a_i^2 + a_j^2 - \text{\rm tr}(A) \, (a_i+a_j) \\ 
&+ \frac{1}{\sigma^2} \, \text{\rm tr}(A)^2 - \frac{2-\sigma}{\sigma} \, |A|^2 \geq 0 
\end{align*}
for $i \neq j$, where $a_i$ denotes the $i$-th eigenvalue of $A$. 

We now verify this inequality. Since $A \geq 0$, we have 
\[\text{\rm tr}(A)^2-|A|^2 = \sum_{p \neq q} a_p a_q \geq 0\] 
and 
\begin{align*} 
&(\text{\rm tr}(A)-a_i-a_j)^2 - (|A|^2-a_i^2-a_j^2) \\ 
&= \sum_{p,q \in \{1,\hdots,n\} \setminus \{i,j\}, \, p \neq q} a_p a_q \geq 0
\end{align*} 
for $i \neq j$. At this point, we distinguish two cases: 

\textit{Case 1:} Suppose first that $\sigma \in (0,\frac{4}{3}]$. In this case, we have 
\begin{align*} 
&(n-2) \, a_i a_j + a_i^2 + a_j^2 - \text{\rm tr}(A) \, (a_i+a_j) \\ 
&+ \frac{1}{\sigma^2} \, \text{\rm tr}(A)^2 - \frac{2-\sigma}{\sigma} \, |A|^2 \\ 
&= (n-3) \, a_i a_j + \frac{(1-\sigma)^2}{\sigma^2} \, \text{\rm tr}(A)^2 \\ 
&+ \frac{4-3\sigma}{2\sigma} \, (\text{\rm tr}(A)^2-|A|^2) \\ 
&+ \frac{1}{2} \, \big ( (\text{\rm tr}(A)-a_i-a_j)^2 - (|A|^2-a_i^2-a_j^2) \big ) \\ 
&\geq 0, 
\end{align*} 
as desired. 

\textit{Case 2:} Suppose now that $\sigma \in [\frac{4}{3},2]$. In this case, we obtain 
\begin{align*}
&(n-2) \, a_i a_j + a_i^2 + a_j^2 - \text{\rm tr}(A) \, (a_i+a_j) \\ 
&+ \frac{1}{\sigma^2} \, \text{\rm tr}(A)^2 - \frac{2-\sigma}{\sigma} \, |A|^2 \\ 
&= \Big ( n-4 + 2 \, \frac{2-\sigma}{\sigma} \Big ) \, a_i a_j + \frac{(2-\sigma)^2}{4\sigma^2} \, \text{\rm tr}(A)^2 \\ 
&+ \frac{3\sigma-4}{\sigma} \, \Big ( \frac{1}{2} \, \text{\rm tr}(A) - a_i-a_j \Big )^2 \\ 
&+ \frac{2-\sigma}{\sigma} \, \big ( (\text{\rm tr}(A)-a_i-a_j)^2 - (|A|^2-a_i^2-a_j^2) \big ) \\ 
&\geq 0. 
\end{align*}
This proves the claim. \\

\textit{Step 3:} In the next step, we show that the condition $\text{\rm tr}(H) \, \text{\rm id} - (n-2\sigma) \, H \geq 0$ is preserved. Using the fact that $\text{\rm Ric}_0(S) = 0$, we compute 
\begin{align*}
&\frac{d}{dt} (\text{\rm tr}(H) \, \text{\rm id} - (n-2\sigma) \, H) \\ 
&= 2 \, S*(\text{\rm tr}(H) \, \text{\rm id} - (n-2\sigma) \, H) + \frac{2}{n} \, \text{\rm scal}(S) \, (\text{\rm tr}(H) \, \text{\rm id} - (n-2\sigma) \, H) \\ 
&+ 4 \, \text{\rm tr}(H) \, (\text{\rm tr}(H) \, \text{\rm id} - (n-2\sigma) \, H) + 4 \, (n-2\sigma) \, H^2 - \frac{4}{n-2\sigma} \, \text{\rm tr}(H)^2 \, \text{\rm id}. 
\end{align*} 
This implies that the condition $\text{\rm tr}(H) \, \text{\rm id} - (n-2\sigma) \, H \geq 0$ is preserved. \\

\textit{Step 4:} We next show that the condition $\text{\rm tr}(H) - \theta \, \text{\rm scal}(S) \geq 0$ is preserved if $\theta$ is sufficiently small. This is trivial if $\theta=0$, so we will only consider the case $\theta>0$. We compute 
\begin{align*} 
\frac{d}{dt} \text{\rm scal}(S) 
&= \frac{2}{n} \, \text{\rm scal}(S)^2 - 2n \, (\text{\rm tr}(H)^2-|H|^2) \\ 
&+ \frac{4n(n-1)}{\sigma(n-2\sigma)} \, \text{\rm tr}(H)^2 - \frac{2n(n-1)(2-\sigma)}{\sigma} \, |H|^2 
\end{align*}
and 
\begin{align*} 
\frac{d}{dt} \text{\rm tr}(H) 
&= \frac{4}{n} \, \text{\rm scal}(S) \, \text{\rm tr}(H) + 4 \, (\text{\rm tr}(H)^2 - |H|^2) \\ 
&- \frac{2n}{\sigma(n-2\sigma)} \, \text{\rm tr}(H)^2 + \frac{2n}{\sigma} \, |H|^2. 
\end{align*}
Hence, if $\text{\rm tr}(H) - \theta \, \text{\rm scal}(S) = 0$, then we obtain 
\begin{align*} 
&\frac{d}{dt} (\text{\rm tr}(H) - \theta \, \text{\rm scal}(S)) \\ 
&= \frac{4}{n} \, \text{\rm scal}(S) \, \text{\rm tr}(H) - \frac{2}{n} \, \theta \, \text{\rm scal}(S)^2 + (4+2n\theta) \, (\text{\rm tr}(H)^2 - |H|^2) \\ 
&- \frac{2n}{\sigma(n-2\sigma)} \, \text{\rm tr}(H)^2 - \frac{4n(n-1)}{\sigma(n-2\sigma)} \, \theta \, \text{\rm tr}(H)^2 \\ 
&+ \frac{2n}{\sigma} \, |H|^2 + \frac{2n(n-1)(2-\sigma)}{\sigma} \, \theta \, |H|^2 \\ 
&= \frac{2}{n\theta} \, \text{\rm tr}(H)^2 + (4+2n\theta) \, (\text{\rm tr}(H)^2 - |H|^2) \\ 
&+ \frac{2+2(n-1)(2-\sigma)\theta}{\sigma} \, (n \, |H|^2-\text{\rm tr}(H)^2) \\ 
&- \frac{4 + 2(n-1)(n+4-2\sigma)\theta}{n-2\sigma} \, \text{\rm tr}(H)^2. 
\end{align*}
If $\theta>0$ is sufficiently small, then the right hand side is nonnegative. Hence, if $\theta$ is sufficiently small, then the condition $\text{\rm tr}(H) - \theta \, \text{\rm scal}(S) \geq 0$ is preserved. \\

\textit{Step 5:} We now show that the sum $S + H \owedge \text{\rm id}$ evolves by the Hamilton ODE. Using the fact that $\text{\rm Ric}_0(S) = 0$, we obtain $\text{\rm Ric}(S) \owedge H = \frac{1}{n} \, \text{\rm scal}(S) \, H \owedge \text{\rm id}$. This gives 
\begin{align*} 
&\frac{d}{dt} \big ( S + H \owedge \text{\rm id} \big ) \\ 
&= Q(S) + (n-2) \, H \owedge H - 2 \, \text{\rm tr}(H) \, H \owedge \text{\rm id} + 2 \, H^2 \owedge \text{\rm id} \\ 
&+ \frac{2}{\sigma(n-2\sigma)} \, \text{\rm tr}(H)^2 \, \text{\rm id} \owedge \text{\rm id} - \frac{2-\sigma}{\sigma} \, |H|^2 \, \text{\rm id} \owedge \text{\rm id} \\
&+ 2 \, (S*H) \owedge \text{\rm id} + \frac{2}{n} \, \text{\rm scal}(S) \, H \owedge \text{\rm id} + 4 \, \text{\rm tr}(H) \, H \owedge \text{\rm id} - 4 \, H^2 \owedge \text{\rm id} \\ 
&- \frac{2}{\sigma(n-2\sigma)} \, \text{\rm tr}(H)^2 \, \text{\rm id} \owedge \text{\rm id} + \frac{2}{\sigma} \, |H|^2 \, \text{\rm id} \owedge \text{\rm id} \\  
&= Q(S) + 2 \, \text{\rm Ric}(S) \owedge H + 2 \, (S*H) \owedge \text{\rm id} \\ 
&+ (n-2) \, H \owedge H + 2 \, \text{\rm tr}(H) \, H \owedge \text{\rm id} \\ 
&- 2 \, H^2 \owedge \text{\rm id} + |H|^2 \, \text{\rm id} \owedge \text{\rm id} \\ 
&= Q(S + H \owedge \text{\rm id}) 
\end{align*} 
in view of Lemma \ref{algebra}. This shows that the cone $\mathcal{C}_{\sigma,\theta}$ is invariant under the Hamilton ODE for each $\sigma \in (0,2]$ and $\theta \in [0,\bar{\theta}]$. \\

\textit{Step 6:} Finally, to prove the transversality statement, we assume that $\sigma \in (0,1) \cup (1,2)$ and $\theta \in (0,\bar{\theta})$. If $R = S + H \owedge \text{\rm id} \in \mathcal{C}_{\sigma,\theta} \setminus \{0\}$, then $\text{\rm tr}(H) > 0$, hence $\text{\rm tr}(A) > 0$. The argument in Step 2 now implies that the term $T$ has strictly positive curvature operator. Moreover, the calculation in Step 4 gives $\frac{d}{dt} (\text{\rm tr}(H) - \theta \, \text{\rm scal}(S)) > 0$ when $\text{\rm tr}(H) - \theta \, \text{\rm scal}(S)=0$. Consequently, the cone $\mathcal{C}_{\sigma,\theta}$ is transversally invariant away from $0$. This completes the proof. \\

As a consequence of Theorem \ref{pinching.cones}, we obtain a higher-dimensional version of the Hamilton-Ivey pinching estimate for three-dimensional Ricci flow (cf. \cite{Hamilton4},\cite{Ivey}). We begin with an auxiliary result, which is analogous to Theorem 4.1 in \cite{Bohm-Wilking} and Proposition 16 in \cite{Brendle-Schoen}:

\begin{lemma}
\label{auxiliary.result}
Fix a compact interval $[\alpha,\beta] \subset (1,2)$ and $\theta \in (0,\bar{\theta})$. Assume that $F_0$ is a closed set which is invariant under the Hamilton ODE $\frac{d}{dt} R = Q(R)$ and satisfies 
\[F_0 \subset \{R: R + h \, \text{\rm id} \owedge \text{\rm id} \in \mathcal{C}_{\sigma,\theta}\}\] 
for some $\sigma \in [\alpha,\beta]$ and some $h>0$. Then there exists a positive real number $\varepsilon$, depending only on $\alpha$, $\beta$, $\theta$, and $n$, such that the set 
\[F_1 = \{R \in F_0: R + 2h \, \text{\rm id} \owedge \text{\rm id} \in \mathcal{C}_{\sigma-\varepsilon,\theta}\}\] 
is invariant under the Hamilton ODE $\frac{d}{dt} R = Q(R)$. 
\end{lemma}

\textbf{Proof.} 
Since the cones $\mathcal{C}_{\sigma,\theta}$ are transversally invariant away from the origin, we can find a large constant $N$, depending only on $\alpha$, $\beta$, $\theta$, and $n$, with the following property: if $R \in \partial \mathcal{C}_{\sigma,\theta}$ for some $\sigma \in [\alpha-\frac{1}{N},\beta]$ and $\text{\rm scal}(R) > N$, then $Q(R-2 \, \text{\rm id} \owedge \text{\rm id})$ lies in the interior of the tangent cone to $\mathcal{C}_{\sigma,\theta}$ at $R$. Moreover, we can find a real number $\varepsilon \in (0,\frac{1}{N})$, depending only on $\alpha$, $\beta$,  $\theta$, and $n$, such that 
\[\{R: R + \text{\rm id} \owedge \text{\rm id} \in \mathcal{C}_{\sigma,\theta}\} \cap \{R: \text{\rm scal}(R) \leq N\} \subset \{R: R + 2 \, \text{\rm id} \owedge \text{\rm id} \in \mathcal{C}_{\sigma-\varepsilon,\theta}\}\] 
for all $\sigma \in [\alpha,\beta]$. 

We claim that $\varepsilon$ has the desired property. To prove this, suppose that $F_0$ is a set which is invariant under the Hamilton ODE and satisfies 
\[F_0 \subset \{R: R + h \, \text{\rm id} \owedge \text{\rm id} \in \mathcal{C}_{\sigma,\theta}\}\] 
for some $\sigma \in [\alpha,\beta]$ and some $h>0$. We claim that the set  
\[F_1 = \{R \in F_0: R + 2h \, \text{\rm id} \owedge \text{\rm id} \in \mathcal{C}_{\sigma-\varepsilon,\theta}\}.\] 
is invariant under the Hamilton ODE. To see this, let $R(t)$ be a solution of the ODE such that $R(0) \in F_1$. Clearly, $R(t) \in F_0$ for all $t \geq 0$. This gives $R(t)+h \, \text{\rm id} \owedge \text{\rm id} \in \mathcal{C}_{\sigma,\theta}$ for all $t \geq 0$. We claim that $R(t) + 2h \, \text{\rm id} \owedge \text{\rm id} \in \mathcal{C}_{\sigma-\varepsilon,\theta}$ for all $t \geq 0$. Suppose this is false. Let $\tau = \inf \{t \geq 0: R(t) + 2h \, \text{\rm id} \owedge \text{\rm id} \notin \mathcal{C}_{\sigma-\varepsilon,\theta}\}$. Clearly, $R(\tau)+2h \, \text{\rm id} \owedge \text{\rm id} \in \partial \mathcal{C}_{\sigma-\varepsilon,\theta}$. In view of our choice of $\varepsilon$, we have $\text{\rm scal}(R(\tau)) \geq Nh$, hence $\text{\rm scal}(R(\tau)+2h \, \text{\rm id} \owedge \text{\rm id})>Nh$. In view of our choice of $N$, $Q(R(\tau))$ lies in the interior of the tangent cone to $\mathcal{C}_{\sigma-\varepsilon,\theta}$ at $R(\tau)+2h \, \text{\rm id} \owedge \text{\rm id}$. This contradicts the definition of $\tau$. Thus, $R(t) + 2h \, \text{\rm id} \owedge \text{\rm id} \in \mathcal{C}_{\sigma-\varepsilon,\theta}$ for all $t \geq 0$. This gives $R(t) \in F_1$ for all $t \geq 0$. \\

\begin{theorem}
\label{pinching.estimate}
Let us fix real numbers $\sigma_0 \in (1,2)$ and $\theta \in (0,\bar{\theta})$. Then there exists a concave and increasing function $f: [0,\infty) \to [0,\infty)$ such that $f(s)=\frac{s}{n-2\sigma_0}$ for $s$ sufficiently small, $\lim_{s \to \infty} \frac{f(s)}{s} = \frac{1}{n-2}$, and the set 
\begin{align*} 
&\{R = S + H \owedge \text{\rm id}: S \in PIC2, \, \text{\rm Ric}_0(S)=0, \\ 
&\hspace{32mm} f(\text{\rm tr}(H)) \, \text{\rm id} - H \geq 0, \\ 
&\hspace{32mm} \text{\rm tr}(H) - \theta \, \text{\rm scal}(S) \geq 0\} 
\end{align*}
is preserved by the Hamilton ODE.
\end{theorem}

\textbf{Proof.} 
Using Lemma \ref{auxiliary.result}, we can construct a decreasing sequence $\sigma_0 > \sigma_1 > \sigma_2 > \hdots$ with $\lim_{j \to \infty} \sigma_j = 1$ and a sequence of invariant sets $F_j$ such that $F_0 = \mathcal{C}_{\sigma_0,\theta}$ and 
\[F_j = \{R \in F_{j-1}: R + 2^{j-1} \, \text{\rm id} \owedge \text{\rm id} \in \mathcal{C}_{\sigma_j,\theta}\}\] 
for $j \geq 1$. Therefore, the intersection 
\[F = \bigcap_{j \in \mathbb{N}} F_j = \mathcal{C}_{\sigma_0,\theta} \cap \bigcap_{j \in \mathbb{N}} \{R: R + 2^{j-1} \, \text{\rm id} \owedge \text{\rm id} \in \mathcal{C}_{\sigma_j,\theta}\}\] 
is invariant under the Hamilton ODE. On the other hand, we may write 
\begin{align*} 
F = &\{R = S + H \owedge \text{\rm id}: S \in PIC2, \, \text{\rm Ric}_0(S)=0, \\ 
&\hspace{32mm} f(\text{\rm tr}(H)) \, \text{\rm id} - H \geq 0, \\ 
&\hspace{32mm} \text{\rm tr}(H) - \theta \, \text{\rm scal}(S) \geq 0\}, 
\end{align*} 
where $f(s) := \min \big \{ \frac{s}{n-2\sigma_0},\inf_{j \in \mathbb{N}} \frac{s+2^j \sigma_j}{n-2\sigma_j} \big \}$. Clearly, $f$ is concave, $f(s)=\frac{s}{n-2\sigma_0}$ for $s$ sufficiently small, and $\lim_{s \to \infty} \frac{f(s)}{s} = \frac{1}{n-2}$. This proves the assertion. \\

\section{Splitting Theorems}

\label{splitting.theorems}

In this section, we collect two splitting theorems, which will play a key role in the subsequent arguments. The following result is a higher-dimensional analogue of Hamilton's crucial splitting theorem (cf. \cite{Hamilton5}, Theorem C5.1):

\begin{proposition}
\label{splitting.1}
Let $(M,g(t))$, $t \in (0,T)$, be a (possibly incomplete) solution to the Ricci flow whose curvature tensor lies in the cone $\mathcal{C}_{1,\theta}$ for some $\theta \in (0,\bar{\theta})$. Moreover, suppose that there exists a point $(p_0,t_0)$ in space-time with the property that the curvature tensor lies on the boundary of the PIC2 cone. Then the manifold $(M,g(t_0))$ is either flat or it is locally isometric to a subset of the round cylinder $S^{n-1} \times \mathbb{R}$.
\end{proposition}

\textbf{Proof.} 
By assumption, the curvature tensor $R$ lies in the cone $\mathcal{C}_{1,\theta}$. This is equivalent to saying that 
\[R - \frac{1}{n-2} \, \text{\rm Ric}_0 \owedge \text{\rm id} - \frac{1}{n} \, \frac{\theta}{1+2(n-1)\theta} \, \text{\rm scal} \, \text{\rm id} \owedge \text{\rm id} \in PIC2\] 
and 
\[|v|^2 \, R - \frac{1}{n-2} \, |v|^2 \, \text{\rm Ric}_0 \owedge \text{\rm id} - \frac{1}{2} \, \text{\rm Ric}_0(v,v) \, \text{\rm id} \owedge \text{\rm id} \in PIC2\] 
for every tangent vector $v$. Let $E$ denote the total space of the vector bundle $TM \oplus T^{\mathbb{C}} M \oplus T^{\mathbb{C}} M$ over $M \times (0,T)$. Moreover, let $\Omega$ be the set of all triplets $(v,z,w) \in E$ such that $v \neq 0$ and $z,w$ are linearly independent. We define a function $\varphi: E \to \mathbb{R}$ by 
\begin{align*} 
\varphi(v,z,w) 
&:= |v|^2 \, R(z,w,\bar{z},\bar{w}) - \frac{1}{n-2} \, |v|^2 \, (\text{\rm Ric}_0 \owedge \text{\rm id})(z,w,\bar{z},\bar{w}) \\ 
&- \frac{1}{2} \, \text{\rm Ric}_0(v,v) \, (\text{\rm id} \owedge \text{\rm id})(z,w,\bar{z},\bar{w}). 
\end{align*}
Clearly, $\varphi$ is nonnegative since $R \in \mathcal{C}_{1,\theta}$. Moreover, $\varphi$ satisfies an inequality of the form 
\[D_t \varphi \geq \mathcal{L} \varphi + L \, \inf_{|\xi| \leq 1} (D^2 \varphi)(\xi,\xi) - L \, |d\varphi| - L \, \varphi,\] 
on $\Omega$, where $\mathcal{L}$ denotes the horizontal Laplacian on $E$ and $L$ is a positive function in $L_{loc}^\infty(\Omega)$. Using Bony's strict maximum principle for degenerate elliptic equations, we conclude that the set $\{(v,z,w) \in \Omega: \varphi(v,z,w)=0\}$ is invariant under parallel transport, for every fixed time (see \cite{Bony} and \cite{Brendle-book}, Section 9).   

We now consider the given time $t_0$, and define 
\begin{align*} 
S &:= \{v \in TM: |v|^2 \, R - \frac{1}{n-2} \, |v|^2 \, \text{\rm Ric}_0 \owedge \text{\rm id} - \frac{1}{2} \, \text{\rm Ric}_0(v,v) \, \text{\rm id} \owedge \text{\rm id} \in \partial PIC2\} \\ 
&= \{v \in TM: \text{\rm there exist linearly independent} \\ 
&\hspace{25mm} \text{\rm vectors $z,w \in T^{\mathbb{C}} M$ such that $\varphi(v,z,w) = 0$}\}. 
\end{align*}
In view of the discussion above, this set is invariant under parallel transport. Moreover, it is easy to see that each fiber of $S$ is a linear subspace of the tangent space of $(M,g(t_0))$. Indeed, the fiber of $S$ over any given point either equals $\{0\}$, or it equals the eigenspace corresponding to the largest eigenvalue of the Ricci tensor at that point. Consequently, $S$ defines a parallel subbundle of $TM$. There are four possibilities:

\textit{Case 1:} The parallel subbundle $S$ has rank $n$. In this case, we have $\text{\rm Ric}_0=0$ and $R \in \partial PIC2$ at each point on $(M,g(t_0))$. Since $R - \frac{1}{n-2} \, \text{\rm Ric}_0 \owedge \text{\rm id} - \frac{1}{n} \, \frac{\theta}{1+2(n-1)\theta} \, \text{\rm scal} \, \text{\rm id} \owedge \text{\rm id} \in PIC2$, we conclude that the scalar curvature of $(M,g(t_0))$ vanishes. Thus, $(M,g(t_0))$ is flat.

\textit{Case 2:} The parallel subbundle $S$ has rank $n-1$. In this case, $(M,g(t_0))$ locally splits as a product of two manifolds of dimension $n-1$ and $1$. By Proposition \ref{C1}, $(M,g(t_0))$ is locally isometric to a subset of the round cylinder $S^{n-1} \times \mathbb{R}$. 

\textit{Case 3:} The parallel subbundle $S$ has rank $k \in \{1,\hdots,n-2\}$. In this case, $(M,g(t_0))$ locally splits as a product of two manifolds of dimension $k$ and $n-k$. Let $\{e_1,\hdots,e_n\}$ be an orthonormal basis of the tangent space at $(p_0,t_0)$ with the property that the fiber of $S$ is spanned by $\{e_1,\hdots,e_k\}$. Clearly, $R(e_1,e_n,e_1,e_n) = 0$ since $(M,g(t_0))$ locally splits as a product. This implies 
\begin{align*} 
0 &\leq \varphi(e_1,e_1,e_n) \\ 
&= R(e_1,e_n,e_1,e_n) - \frac{1}{n-2} \, (\text{\rm Ric}_0(e_1,e_1)+\text{\rm Ric}_0(e_n,e_n)) - \text{\rm Ric}_0(e_1,e_1) \\ 
&= -\frac{1}{n-2} \, ((n-1) \, \text{\rm Ric}_0(e_1,e_1)+\text{\rm Ric}_0(e_n,e_n)). 
\end{align*} 
On the other hand, we know that the eigenspace corresponding to the largest eigenvalue of the Ricci tensor is spanned by $\{e_1,\hdots,e_k\}$. Since $k \leq n-2$, we obtain 
\[(n-1) \, \text{\rm Ric}_0(e_1,e_1) + \text{\rm Ric}_0(e_n,e_n) > \sum_{i=1}^{n-1} \text{\rm Ric}_0(e_i,e_i)+\text{\rm Ric}_0(e_n,e_n) = 0.\] 
This is a contradiction. 

\textit{Case 4:} The parallel subbundle $S$ has rank $0$. In this case, 
\[|v|^2 \, R - \frac{1}{n-2} \, |v|^2 \, \text{\rm Ric}_0 \owedge \text{\rm id} - \frac{1}{2} \, \text{\rm Ric}_0(v,v) \, \text{\rm id} \owedge \text{\rm id}\] 
lies in the interior of the PIC2 cone for each $v \neq 0$. If we choose $v$ to be an eigenvector corresponding to the largest eigenvalue of the Ricci tensor, then $\frac{1}{n-2} \, |v|^2 \, \text{\rm Ric}_0 + \frac{1}{2} \, \text{\rm Ric}_0(v,v) \, \text{\rm id}$ is weakly two-positive, hence $\frac{1}{n-2} \, |v|^2 \, \text{\rm Ric}_0 \owedge \text{\rm id} + \frac{1}{2} \, \text{\rm Ric}_0(v,v) \, \text{\rm id} \owedge \text{\rm id} \in PIC2$. Putting these facts together, we conclude that $R$ lies in the interior of the PIC2 cone, contrary to our assumption. This completes the proof. \\

The next result is an adaptation of a result of Perelman \cite{Perelman1} (see also \cite{Chen-Zhu}, Lemma 3.1).

\begin{proposition}
\label{splitting.2}
Let $(M,g)$ be a complete noncompact manifold whose curvature tensor lies in the cone $\mathcal{C}_{1,\theta}$ and in the interior of the PIC2 cone. Let us fix a point $p \in M$ and let $p_j$ be a sequence of points such that $d(p,p_j) \to \infty$. Moreover, let $\lambda_j$ be a sequence of positive real numbers such that $\lambda_j \, d(p,p_j)^2 \to \infty$. If the rescaled manifolds $(M,\lambda_j g,p_j)$ converge in the Cheeger-Gromov sense to a smooth, non-flat limit $Y$, then $Y$ is isometric to a round cylinder $S^{n-1} \times \mathbb{R}$.
\end{proposition} 

\textbf{Proof.} 
By Proposition 2.3 in \cite{Chen-Zhu}, the limit $Y$ splits as a product $Y = X \times \mathbb{R}$, where $X$ is a non-flat $(n-1)$-dimensional manifold. Since the curvature of $Y$ lies in the cone $\mathcal{C}_{1,\theta}$, we conclude that $X$ has constant curvature by Proposition \ref{C1}. Consequently, $X$ is isometric to a space form $S^{n-1}/\Gamma$. A result due to Hamilton (see Theorem \ref{incompressibility.noncompact.case}) then implies that the cross-section $X$ is incompressible in $M$. Since $M$ is diffeomorphic to $\mathbb{R}^n$ by the soul theorem (cf. \cite{Cheeger-Gromoll2}), it follows that $\Gamma$ is trivial. Thus $X$ is isometric to a round sphere $S^{n-1}$, and $Y$ is isometric to a round cylinder $S^{n-1} \times \mathbb{R}$. \\

For later use, we also recall the following result due to Perelman: 

\begin{proposition}[cf. G.~Perelman \cite{Perelman1}]
\label{nonexistence.of.small.necks}
Let $(M,g)$ be a complete Riemannian manifold whose curvature tensor lies in the cone $\mathcal{C}_{1,\theta}$ and in the interior of the PIC2 cone. Moreover, suppose that $(M,g)$ is $\kappa$-noncollapsed, and that the covariant derivatives of the Riemann curvature tensor satisfy the pointwise estimates $|D^m R| \leq \eta \, \text{\rm scal}^{\frac{m}{2}+1}$, $1 \leq m \leq 8$, at all points where the scalar curvature is sufficiently large. Then $(M,g)$ has bounded curvature. 
\end{proposition}

\textbf{Proof.} 
Suppose that $(M,g)$ does not have bounded curvature. Using a standard point-picking lemma, we can find a sequence of points $x_j$ such that $Q_j := \text{\rm scal}(x_j) \geq j$ and 
\[\sup_{x \in B(x_j,j \, Q_j^{-\frac{1}{2}})} \text{\rm scal}(x) \leq 4 Q_j.\] 
We now dilate the manifold $(M,g)$ around the point $x_j$ by the factor $Q_j$. Using the noncollapsing assumption and the curvature derivative estimates, we are able to take a limit in $C_{loc}^6$. The limit manifold $(M^\infty,g^\infty)$ is complete, non-flat, and has bounded curvature. By Proposition 2.3 in \cite{Chen-Zhu}, the limit $Y$ splits as a product $Y = X \times \mathbb{R}$, where $X$ is a non-flat $(n-1)$-dimensional manifold. As above, Proposition \ref{C1} implies that $X$ is isometric to a space form, and Theorem \ref{incompressibility.noncompact.case} gives that $X$ is a round sphere. Consequently, the original manifold $(M,g)$ contains a sequence of necks with radii converging to $0$. But this is impossible in a manifold with nonnegative sectional curvature (see \cite{Chen-Zhu}, Proposition 2.2). \\

\section{Ancient $\kappa$-solutions with $\theta$-pinched curvature}

\label{ancient.solutions}

In this section, we discuss how the arguments in Section 11 of Perelman's first paper can be extended to higher dimensions. Through this section, we fix an integer $n \geq 5$ and an arbitrary positive constant $\theta \in (0,\bar{\theta})$. Our goal is to analyze ancient solutions whose curvature tensor lies in the cone $\mathcal{C}_{1,\theta}$. We will use the following terminology:

\begin{definition}
An ancient $\kappa$-solution with $\theta$-pinched curvature is a non-flat ancient solution to the Ricci flow of dimension $n$ which is complete with bounded curvature; is $\kappa$-noncollapsed on all scales; and has curvature in the cone $\mathcal{C}_{1,\theta}$. 
\end{definition}

In view of \cite{Brendle2}, Hamilton's trace Harnack inequality holds in this setting: 

\begin{theorem}[cf. R.~Hamilton \cite{Hamilton3}] 
\label{harnack}
Let $(M,g(t))$ be an ancient $\kappa$-solution with $\theta$-pinched curvature. Then 
\[\frac{\partial}{\partial t} \text{\rm scal} + 2 \, \langle \nabla \text{\rm scal},v \rangle + 2 \, \text{\rm Ric}(v,v) \geq 0\] 
for every tangent vector $v$. In particular, the scalar curvature is monotone increasing at each point. 
\end{theorem} 

\textbf{Proof.} 
This was established in a seminal paper of Hamilton \cite{Hamilton3} under the stronger assumption of nonnegative curvature operator. In \cite{Brendle2}, we showed that Hamilton's Harnack estimate holds under the weaker assumption that the curvature tensor lies in the cone PIC2. Since the cone $\mathcal{C}_{1,\theta}$ is contained in the cone PIC2, the assertion follows. \\

Integrating the trace Harnack inequality along paths in space-time gives the following result:

\begin{corollary}
\label{integrated.harnack}
Let $(M,g(t))$ be an ancient $\kappa$-solution with $\theta$-pinched curvature. Then
\[\text{\rm scal}(x_1,t_1) \leq \exp \Big ( \frac{d_{g(t_1)}(x_1,x_2)^2}{2(t_2-t_1)} \Big ) \, \text{\rm scal}(x_2,t_2)\] 
whenever $t_1 < t_2$.
\end{corollary}

We next recall a key result from Perelman's first paper:

\begin{theorem}[cf. G.~Perelman \cite{Perelman1}, Corollary 11.6]
\label{perelman.cor.11.6}
For every $w>0$, there exist positive constants $B$ and $C$ with the following property: Suppose that $(M,g(t))$, $t \in [0,T]$, is a solution to the Ricci flow so that the ball $B_{g(T)}(x_0,r_0)$ is compactly contained in $M$. Moreover, suppose that for each $t \in [0,T]$, the curvature tensor of $g(t)$ lies in the PIC2 cone, and $\text{\rm vol}_{g(t)}(B_{g(t)}(x_0,r_0)) \geq wr_0^n$. Then $\text{\rm scal}(x,t) \leq C r_0^{-2} + Bt^{-1}$ for all $t \in (0,T]$ and all $x \in B_{g(t)}(x_0,\frac{1}{4} \, r_0)$.
\end{theorem}

\textbf{Proof.} 
The only difference to the statement in Perelman's paper is that we have replaced the assumption that $g(t)$ has nonnegative curvature operator by the weaker PIC2 assumption. This does not affect the proof. \\

In the next step, we establish an analogue of Perelman's longrange curvature estimate for ancient $\kappa$-solutions in dimension $3$. Perelman's estimate was adapted to the four-dimensional case in \cite{Chen-Zhu}, Proposition 3.3.

\begin{theorem}[cf. G.~Perelman \cite{Perelman1}, Section 11.7]
\label{longrange.curvature.estimate}
Given $\kappa>0$, there exists a large positive constant $\eta$ and a positive function $\omega: [0,\infty) \to (0,\infty)$ (depending on $\kappa$) such that the following holds: Let $(M,g(t))$ be an ancient $\kappa$-solution with $\theta$-pinched curvature. Then 
\[\text{\rm scal}(x,t) \leq \text{\rm scal}(y,t) \, \omega(\text{\rm scal}(y,t) \, d_{g(t)}(x,y)^2)\] 
for all points $x,y \in M$ and all $t$. Furthermore, we have the pointwise estimates $|D^m R| \leq \eta \, \text{\rm scal}^{\frac{m}{2}+1}$ for $1 \leq m \leq 8$, where $R$ denotes the Riemann curvature tensor.
\end{theorem}

\textbf{Proof.} 
We sketch the argument, following Section 11.7 in Perelman's paper \cite{Perelman1} (see also \cite{Chen-Zhu}, Proposition 3.3). The second statement follows immediately from the first statement together with Shi's interior derivative estimate. Thus, it suffices to prove the first statement. Without loss of generality, we may assume that $t=0$. Moreover, by scaling we may assume that $\text{\rm scal}(y,0) = 1$. Let $A$ denote the set of all points $x \in M$ such that $\text{\rm scal}(x,0)+1 \geq d_{g(0)}(y,x)^{-2}$. Moreover, let $z \in A$ be a point which has minimal distance from $y$ with respect to the metric $g(0)$ among all points in $A$. Clearly, $\text{\rm scal}(z,0)+1 = d_{g(0)}(y,z)^{-2}$.

Let $p$ denote the mid-point of the minimizing geodesic in $(M,g(0))$ joining $y$ and $z$. It follows from the definition of $z$ that $B_{g(0)}(p,\frac{1}{4} \, d_{g(0)}(y,z)) \cap A = \emptyset$. In other words, $\text{\rm scal}(x,0)+1 \leq d_{g(0)}(y,x)^{-2}$ for all $x \in B_{g(0)}(p,\frac{1}{4} \, d_{g(0)}(y,z))$. Using the Harnack inequality, we obtain  
\begin{align*} 
\sup_{x \in B_{g(t)}(p,\frac{1}{4} \, d_{g(0)}(y,z))} \text{\rm scal}(x,t) 
&\leq \sup_{x \in B_{g(0)}(p,\frac{1}{4} \, d_{g(0)}(y,z))} \text{\rm scal}(x,0) \\ 
&\leq 16 \, d_{g(0)}(y,z)^{-2} 
\end{align*} 
for all $t \in (-\infty,0]$. The noncollapsing property gives 
\[\text{\rm vol}_{g(t)} \big ( B_{g(t)}(p,\frac{1}{4} \, d_{g(0)}(y,z)) \big ) \geq \kappa \, (\frac{1}{4} \, d_{g(0)}(y,z))^n\] 
for all $t \in (-\infty,0]$. Using Theorem \ref{perelman.cor.11.6}, we conclude that there exists a positive and increasing function $\omega: [0,\infty) \to (0,\infty)$ such that 
\[\sup_{x \in B_{g(0)}(p,r)} \text{\rm scal}(x,0) \leq d_{g(0)}(y,z)^{-2} \, \omega(d_{g(0)}(y,z)^{-1} \, r)\] 
for all $r \geq 0$. The Harnack inequality gives 
\begin{align*} 
\sup_{x \in B_{g(0)}(p,d_{g(0)}(y,z))} \text{\rm scal}(x,t) 
&\leq \sup_{x \in B_{g(0)}(p,d_{g(0)}(y,z))} \text{\rm scal}(x,0) \\ 
&\leq d_{g(0)}(y,z)^{-2} \, \omega(1)
\end{align*}
for all $t \in (-\infty,0]$. In particular, there exists a positive constant $\beta$ such that $d_{g(t)}(y,z) \leq 2 \, d_{g(0)}(y,z)$ for all $t \in [-\beta \, d_{g(0)}(y,z)^2,0]$. Moreover, by choosing $\beta$ sufficiently small, we can arrange that $\text{\rm scal}(z,0) \leq \text{\rm scal}(z,t) + \frac{1}{2} \, d_{g(0)}(y,z)^{-2}$ for all $t \in [-\beta \, d_{g(0)}(y,z)^2,0]$. (This follows from Shi's interior derivative estimate.) Applying the Harnack inequality with $t = -\beta \, d_{g(0)}(y,z)^2$ yields 
\begin{align*} 
\frac{1}{2} \, d_{g(0)}(y,z)^{-2} - 1 
&= \text{\rm scal}(z,0) - \frac{1}{2} \, d_{g(0)}(y,z)^{-2} \\ 
&\leq \text{\rm scal}(z,t) \\ 
&\leq \exp \Big ( -\frac{d_{g(t)}(y,z)^2}{2t} \Big ) \, \text{\rm scal}(y,0) \\ 
&\leq \exp \Big ( -\frac{2 \, d_{g(0)}(y,z)^2}{t} \Big ) \, \text{\rm scal}(y,0) \\ 
&= \exp \Big ( \frac{2}{\beta} \Big ). 
\end{align*} 
This implies $d_{g(0)}(y,z)^{-2} \leq 4 \, e^{\frac{2}{\beta}}$. Moreover, we have $d_{g(0)}(y,p) = \frac{1}{2} \, d_{g(0)}(y,z) \leq \frac{1}{2}$. Putting these facts together, we obtain 
\begin{align*} 
\sup_{x \in B_{g(0)}(y,r)} \text{\rm scal}(x,0) 
&\leq \sup_{x \in B_{g(0)}(p,r+1)} \text{\rm scal}(x,0) \\ 
&\leq d_{g(0)}(y,z)^{-2} \, \omega(d_{g(0)}(y,z)^{-1} \, (r+1)) \\ 
&\leq 4 \, e^{\frac{2}{\beta}} \, \omega(2 \, e^{\frac{1}{\beta}} \, (r+1)) 
\end{align*} 
for all $r \geq 0$. This proves the assertion. \\

As a corollary, we obtain a higher dimensional version of Perelman's compactness theorem for ancient $\kappa$-solutions:

\begin{corollary}[cf. G.~Perelman \cite{Perelman1}, Section 11.7]
\label{compactness.thm.for.ancient.solutions}
Fix $\kappa>0$ and $\theta \in (0,\bar{\theta})$, and let $(M^{(j)},g^{(j)}(t))$ be a sequence of ancient $\kappa$-solutions with $\theta$-pinched curvature. Suppose that each solution is defined for $t \in (-\infty,0]$, and that $\text{\rm scal}(x_j,0) = 1$ for some point $x_j$. Then, after passing to a subsequence if necessary, the sequence $(M^{(j)},g^{(j)}(t),x_j)$ converges in the Cheeger-Gromov sense to an ancient $\kappa$-solution with $\theta$-pinched curvature.
\end{corollary}

\textbf{Proof.} 
It follows from Theorem \ref{longrange.curvature.estimate} that, after passing to a subsequence if necessary, the sequence $(M^{(j)},g^{(j)}(t),x_j)$ converges in the Cheeger-Gromov sense to a smooth ancient solution $(M^\infty,g^\infty(t))$. Clearly, $(M^\infty,g^\infty(t))$ is $\kappa$-noncollapsed and has $\theta$-pinched curvature. Moreover, Theorem \ref{longrange.curvature.estimate} implies that $(M^{(j)},g^{(j)}(t))$ satisfies $|D^m R| \leq \eta \, \text{\rm scal}^{\frac{m}{2}+1}$ for $1 \leq m \leq 8$, where $\eta$ is a positive constant which does not depend on $j$. Hence, these estimates also hold on the limiting ancient solution $(M^\infty,g^\infty(t))$. Proposition \ref{nonexistence.of.small.necks} now implies that the limiting ancient solution has bounded curvature. \\

In the remaining part of this section, we establish a universal noncollapsing property and a Canonical Neighborhood Theorem for ancient $\kappa$-solutions. This was first established by Perelman \cite{Perelman2} in the three-dimensional case, and adapted to dimension $4$ in \cite{Chen-Zhu}. We begin with the noncompact case:

\begin{theorem}[cf. G.~Perelman \cite{Perelman1}, Corollary 11.8; Chen-Zhu \cite{Chen-Zhu}, Proposition 3.4]
\label{canonical.neighborhood.preparation.1}
Given $\varepsilon>0$, there exist positive constants $C_1 = C_1(n,\theta,\varepsilon)$ and $C_2 = C_2(n,\theta,\varepsilon)$ such that the following holds: Let $(M,g(t))$ be a noncompact, non-flat ancient $\kappa$-solution with $\theta$-pinched curvature which is not locally isometric to a round cylinder. Given any point $(x_0,t_0)$ in space-time there exists an open neighborhood $B$ of $x_0$ such that $B_{g(t_0)}(x_0,C_1^{-1} \, \text{\rm scal}(x_0,t_0)^{-\frac{1}{2}}) \subset B \subset B_{g(t_0)}(x_0,C_1 \, \text{\rm scal}(x_0,t_0)^{-\frac{1}{2}})$ and $C_2^{-1} \, \text{\rm scal}(x_0,t_0) \leq \text{\rm scal}(x,t_0) \leq C_2 \, \text{\rm scal}(x_0,t_0)$ for all $x \in B$. Moreover, one of the following statements holds: 
\begin{itemize}
\item $B$ is an $\varepsilon$-neck.
\item $B$ is an $\varepsilon$-cap in the sense that $B$ is diffeomorphic to a ball and the boundary $\partial B$ is a cross-sectional sphere of an $\varepsilon$-neck. 
\end{itemize}
In particular, $(M,g(t_0))$ is $\kappa_0$-noncollapsed for some universal constant $\kappa_0 = \kappa_0(n)$.
\end{theorem}

\textbf{Proof.} 
Without loss of generality, we may assume that $t_0=0$. Moreover, we may assume that $x_0$ does not lie at the center of an $\varepsilon$-neck. Since $(M,g(t))$ is not locally isometric to a round cylinder, it follows from Proposition \ref{splitting.1} that the curvature tensor of $(M,g(t))$ lies in the interior of the PIC2 cone. In particular, $M$ is diffeomorphic to $\mathbb{R}^n$ by the soul theorem. Let $M_\varepsilon$ denote the set of all points in $(M,g(0))$ which do not lie at the center of an $\frac{\varepsilon}{2}$-neck. Clearly, $x_0 \in M_\varepsilon$ since $x_0$ does not lie at the center of an $\varepsilon$-neck.

\textit{Step 1:} We first show that the closure of $M_\varepsilon$ is compact. Suppose this is false. Then there exists a sequence of points $x_j \in M_\varepsilon$ such that $d_{g(0)}(x_0,x_j) \to \infty$. Since $\text{\rm scal}(x_0,0) > 0$, Theorem \ref{longrange.curvature.estimate} implies that $\liminf_{j \to \infty} \lambda_j \, d_{g(0)}(x_0,x_j)^2 = \infty$, where $\lambda_j = \text{\rm scal}(x_j,0)$. Using Corollary \ref{compactness.thm.for.ancient.solutions}, we conclude that the rescaled manifolds $(M,\lambda_j \, g(0),x_j)$ converge in the Cheeger-Gromov sense to a smooth non-flat limit. Since $\lambda_j \, d_{g(0)}(x_0,x_j)^2 \to \infty$, Proposition \ref{splitting.2} implies that the limit is isometric to a cylinder. In particular, $x_j$ lies on an $\frac{\varepsilon}{2}$-neck if $j$ is sufficiently large. This contradicts the fact that $x_j \in M_\varepsilon$. 

Thus, $M_\varepsilon$ has compact closure. In particular, $M_\varepsilon \neq M$. Since $M_\varepsilon \neq \emptyset$, it follows that $\partial M_\varepsilon \neq \emptyset$. 

\textit{Step 2:} Let us consider an arbitrary point $y \in \partial M_\varepsilon$. Clearly, $y$ lies on an $\varepsilon$-neck in $(M,g(0))$. Using the Harnack inequality, we obtain $\text{\rm scal}(x,t) \leq \text{\rm scal}(x,0) \leq 2 \, \text{\rm scal}(y,0)$ for all $x \in B_{g(0)}(y,\text{\rm scal}(y,0)^{-\frac{1}{2}})$ and all $t \leq 0$. Hence, there exists a small constant $\beta=\beta(n)>0$ such that 
\[\text{\rm vol}_{g(t)} \big ( B_{g(t)}(y,\text{\rm scal}(y,0)^{-\frac{1}{2}}) \big ) \geq \beta \, \text{\rm scal}(y,0)^{-\frac{n}{2}}\] 
for all $t \in [-\beta \, \text{\rm scal}(y,0)^{-1},0]$. Using Theorem \ref{perelman.cor.11.6}, we conclude that 
\[\text{\rm scal}(x,0) \leq \text{\rm scal}(y,0) \, \omega(\text{\rm scal}(y,0) \, d_{g(0)}(x,y)^2)\] 
for all $x \in M$, where $\omega$ is a positive function that does not depend on $\kappa$.

\textit{Step 3:} We again consider an arbitrary point $y \in \partial M_\varepsilon$. Recall that $y$ lies on an $\varepsilon$-neck in $(M,g(0))$. Let $\Sigma_y$ denote the leaf in Hamilton's CMC foliation which passes through the point $y$. Since $M$ is diffeomorphic to $\mathbb{R}^n$, there is a unique bounded connected component of $M \setminus \Sigma_y$, and this connected component is diffeomorphic to a ball in view of the solution of the Schoenflies conjecture in dimension $n \neq 4$. Let us denote this connected component by $\Omega_y$. 

We claim that $\text{\rm scal}(y,0) \, \text{\rm diam}_{g(0)}(\Omega_y)^2 \leq C$, where $C$ depends on $n$, $\theta$, and $\varepsilon$, but not on $\kappa$. To prove this, we argue by contradiction. Suppose that $(M^{(j)},g^{(j)}(t))$ is a sequence of noncompact, non-flat ancient $\kappa_j$-solutions with $\theta$-pinched curvature which are not locally isometric to a round cylinder. Moreover, suppose that $y_j$ is a sequence of points such that $y_j \in \partial M_\varepsilon^{(j)}$ and $\text{\rm scal}(y_j,0) \, \text{\rm diam}_{g^{(j)}(0)}(\Omega_{y_j})^2 \to \infty$, where $\Omega_{y_j}$ denotes the region in $(M^{(j)},g^{(j)}(0))$ which is bounded by the CMC sphere passing through $y_j$. We dilate the manifold $(M^{(j)},g^{(j)}(0))$ around the point $y_j$ by the factor $\text{\rm scal}(y_j,0)$. The curvature estimate established in Step 2 implies that, after passing to a subsequence if necessary, the rescaled manifolds converge to a smooth limit. Moreover, since $\text{\rm scal}(y_j,0) \, \text{\rm diam}_{g^{(j)}(0)}(\Omega_{y_j})^2 \to \infty$, the limiting manifold has at least two ends. Consequently, the limiting manifold splits off a line by the Cheeger-Gromoll splitting theorem \cite{Cheeger-Gromoll1}. Since the curvature tensor lies in the cone $\mathcal{C}_{1,\theta}$ at each point, the limit must be isometric to $(S^{n-1}/\Gamma) \times \mathbb{R}$. If $\Gamma$ is non-trivial, then the manifold $M^{(j)}$ contains a quotient neck if $j$ is sufficiently large. As above, Theorem \ref{incompressibility.noncompact.case} implies that the cross section of such a quotient neck is incompressible in $M^{(j)}$, which is impossible since $M^{(j)}$ is diffeomorphic to $\mathbb{R}^n$. Thus, $\Gamma$ is trivial. Thus, the limit is isometric to the round cylinder $S^{n-1} \times \mathbb{R}$. Hence, if $j$ is sufficiently large, then $y_j$ lies in the interior of the set $M_\varepsilon^{(j)}$. This contradicts our choice of $y_j$. 

\textit{Step 4:} Combining the curvature estimate in Step 2 with the diameter estimate in Step 3 gives $\text{\rm scal}(x,0) \leq C \, \text{\rm scal}(y,0)$ for all $y \in \partial M_\varepsilon$ and all $x \in \Omega_y$. Here, $C$ is a positive constant that depends only on $\varepsilon$, but not on $\kappa$. Moreover, the Harnack inequality (cf. Corollary \ref{integrated.harnack} above) implies that $\text{\rm scal}(x,0) \geq \frac{1}{C} \, \text{\rm scal}(y,0)$ for all $y \in \partial M_\varepsilon$ and all $x \in \Omega_y$, where $C$ depends only on $\varepsilon$, but not on $\kappa$. 

\textit{Step 5:} Finally, we observe that the sets $\Omega_y$ are nested. More precisely, given two points $y,y' \in \partial M_\varepsilon$, we either have $\Omega_y \subset \Omega_{y'}$ or $\Omega_{y'} \subset \Omega_y$. Since $\partial M_\varepsilon$ is compact, we can find a point $y_0 \in \partial M_\varepsilon$ such that $\Omega_y \subset \Omega_{y_0}$ for all $y \in \partial M_\varepsilon$. In particular, the set $\partial M_\varepsilon$ is contained in the closure of $\Omega_{y_0}$. Since $M_\varepsilon$ has compact closure, it follows that $M_\varepsilon$ is contained in the closure of $\Omega_{y_0}$. In particular, we have $x_0 \in \Omega_{y_0}$. Since $x_0$ does not lie at the center of an $\varepsilon$-neck, the distance of $x_0$ to the boundary $\partial \Omega_{y_0} = \Sigma_{y_0}$ is bounded from below by $C^{-1} \, \text{\rm scal}(x_0,0)^{-\frac{1}{2}}$. Therefore, $B_{g(0)}(x_0,C^{-1} \, \text{\rm scal}(x_0,0)^{-\frac{1}{2}}) \subset \Omega_{y_0}$. On the other hand, the diameter bound in Step 3 gives $\Omega_{y_0} \subset B_{g(0)}(x_0,C \, \text{\rm scal}(x_0,0)^{-\frac{1}{2}})$. Thus, the set $B := \Omega_{y_0}$ has all the required properties. \\

\begin{theorem}[cf. G.~Perelman \cite{Perelman1}; Chen-Zhu \cite{Chen-Zhu}]
\label{universal.noncollapsing} 
There exists a constant $\kappa_0 = \kappa_0(n,\theta)$ such that the following holds: Let $(M,g(t))$ be an ancient $\kappa$-solution with $\theta$-pinched curvature for some $\kappa>0$. Then either $(M,g(t))$ is $\kappa_0$-noncollapsed for all $t$; or $(M,g(t))$ is a metric quotient of the round sphere $S^n$; or $(M,g(t))$ is a noncompact quotient of the round cylinder $S^{n-1} \times \mathbb{R}$.
\end{theorem}

\textbf{Proof.} 
If $M$ is noncompact, the assertion is a consequence of Theorem \ref{canonical.neighborhood.preparation.1}. Hence, it suffices to consider the case that $M$ is compact. In view of the noncollapsing property, $(M,g(t))$ cannot be a compact quotient of a round cylinder. Hence, Proposition \ref{splitting.1} implies that the curvature tensor of $(M,g(t))$ lies in the interior of the PIC2 cone at each point in space-time. Let $(\bar{M},\bar{g}(t))$ denote the asymptotic shrinking soliton of $(M,g(t))$ (cf. \cite{Perelman1}, Section 11.2). By Perelman's work, $(\bar{M},\bar{g}(t))$ is a non-flat shrinking gradient soliton. By Corollary \ref{compactness.thm.for.ancient.solutions}, $(\bar{M},\bar{g}(t))$ is an ancient $\kappa$-solution with $\theta$-pinched curvature. 

\textit{Case 1:} Suppose that $\bar{M}$ is compact. It is easy to see that a compact shrinking soliton cannot locally split as a product. Consequently, it follows from Proposition \ref{splitting.1} that the curvature tensor of $\bar{M}$ lies in the interior of the PIC2 cone. Using results in \cite{Brendle-Schoen}, we conclude that $\bar{M}$ is isometric to a metric quotient of the round sphere $S^n$. This directly implies that the flow $(M,g(t))$ is a metric quotient of $S^n$. 

\textit{Case 2:} Suppose next that $\bar{M}$ is noncompact. We claim that $\bar{M}$ is noncollapsed with some universal constant. By Theorem \ref{canonical.neighborhood.preparation.1}, the asymptotic shrinking soliton $\bar{M}$ is either $\kappa_0$-noncollapsed for some universal constant $\kappa_0$, or it is isometric to a metric quotient of the round cylinder. Let us examine the latter case. If $\bar{M} = (S^{n-1} \times \mathbb{R})/\Gamma$ and $n$ is odd, there are only finitely many possibilities for the group $\Gamma$, and the resulting quotients are all noncollapsed with a universal constant. On the other hand, if $\bar{M} = (S^{n-1} \times \mathbb{R})/\Gamma$ and $n$ is even, then a result of Hamilton (cf. Theorem \ref{incompressibility.compact.case}) implies the center slice $(S^{n-1} \times \{0\})/\Gamma$ is incompressible in $M$. However, since $n$ is even and the curvature tensor of $M$ lies in the interior of the PIC2 cone, the fundamental group of $M$ has order at most $2$ by Synge's theorem. Again, this leaves only finitely many possibilities for the group $\Gamma$, and the resulting quotients $(S^{n-1} \times \mathbb{R})/\Gamma$ are noncollapsed with a universal constant. To summarize, we have shown that $\bar{M}$ is noncollapsed with a universal constant. Using Perelman's monotonicity formula for the reduced volume, we can deduce that $(M,g(t))$ is noncollapsed with a universal constant. The proof of this is based on work of Perelman \cite{Perelman1}, Section 7.3; for details, see \cite{Chen-Zhu}, pp.~205--208. \\

Using the universal noncollapsing property in Theorem \ref{universal.noncollapsing}, we can extend the Canonical Neighborhood Theorem for ancient $\kappa$-solutions to the compact case:

\begin{corollary}[cf. G.~Perelman \cite{Perelman2}, Section 1.5]
\label{canonical.neighborhood.preparation.2}
Given $\varepsilon>0$, we can find constants $C_1=C_1(n,\theta,\varepsilon)$ and $C_2=C_2(n,\theta,\varepsilon)$ with the following property: Suppose that $(M,g(t))$ is a non-flat ancient $\kappa$-solution with $\theta$-pinched curvature. Then, for each point $(x_0,t_0)$ in space-time there exists a neighborhood $B$ of $x_0$ such that $B_{g(t_0)}(x_0,C_1^{-1} \, \text{\rm scal}(x_0,t_0)^{-\frac{1}{2}}) \subset B \subset B_{g(t_0)}(x_0,C_1 \, \text{\rm scal}(x_0,t_0)^{-\frac{1}{2}})$ and $C_2^{-1} \, \text{\rm scal}(x_0,t_0) \leq \text{\rm scal}(x,t_0) \leq C_2 \, \text{\rm scal}(x_0,t_0)$ for all $x \in B$. Finally, one of the following statements holds: 
\begin{itemize}
\item $B$ is an $\varepsilon$-neck.
\item $B$ is an $\varepsilon$-cap in the sense that $B$ is diffeomorphic to a ball and the boundary $\partial B$ is a cross-sectional sphere of an $\varepsilon$-neck. 
\item $B$ is a closed manifold diffeomorphic to $S^n/\Gamma$.
\item $B$ is an $\varepsilon$-quotient neck of the form $(S^{n-1} \times [-L,L])/\Gamma$.
\end{itemize}
\end{corollary}

\textbf{Proof.} 
It suffices to prove the assertion for $t_0=0$. We argue by contradiction. Suppose that the assertion is false. Then there exists a sequence of non-flat ancient $\kappa_j$-solutions $(M^{(j)},g^{(j)}(t))$ with $\theta$-pinched curvature and a sequence of points $x_j \in M^{(j)}$ with the following property: if $B$ is a neighborhood of $x_j$ such that $B_{g^{(j)}(0)}(x_j,j^{-1} \, \text{\rm scal}(x_j,0)^{-\frac{1}{2}}) \subset B \subset B_{g^{(j)}(0)}(x_j,j \, \text{\rm scal}(x_j,0)^{-\frac{1}{2}})$ and $j^{-1} \, \text{\rm scal}(x_j,0) \leq \text{\rm scal}(x,0) \leq j \, \text{\rm scal}(x_j,0)$ for all $x \in B$, then $B$ is neither an $\varepsilon$-neck; nor an $\varepsilon$-cap; nor a closed manifold diffeomorphic to $S^n/\Gamma$; nor an $\varepsilon$-quotient neck. In particular, $(M^{(j)},g^{(j)}(t))$ cannot be a noncompact quotient of a round cylinder. By Theorem \ref{universal.noncollapsing}, $(M^{(j)},g^{(j)}(t))$ is $\kappa_0$-noncollapsed for some uniform constant $\kappa_0$ which does not depend on $j$.

By scaling, we may assume that $\text{\rm scal}(x_j,0) = 1$ for each $j$. We now apply Corollary \ref{compactness.thm.for.ancient.solutions} to the sequence $(M^{(j)},g^{(j)}(t))$. Hence, after passing to a subsequence, we may assume that the sequence $(M^{(j)},g^{(j)}(t),x_j)$ converges to a non-flat ancient $\kappa_0$-solution with $\theta$-pinched curvature, which we denote by $(M^\infty,g^\infty(t))$. Moreover, the points $x_j \in M^{(j)}$ will converge to a point $x_\infty \in M^\infty$. We distinguish two cases:

\textit{Case 1:} Suppose that $M^\infty$ is compact. In this case, the diameter of $(M^{(j)},g^{(j)}(0))$ is bounded from above by a uniform constant which is independent of $j$. Hence, if $j$ is sufficiently large, then $B^{(j)} := M^{(j)}$ is a neighborhood of the point $x_j$ which satisfies $B_{g^{(j)}(0)}(x_j,j^{-1}) \subset B^{(j)} \subset B_{g^{(j)}(0)}(x_j,j)$ and $j^{-1} \leq \text{\rm scal}(x,0) \leq j$ for all $x \in B^{(j)}$. Furthermore, results in \cite{Brendle-Schoen} imply that $B^{(j)}$ is diffeomorphic to a space form. This contradicts the definition of $x_j$. 

\textit{Case 2:} Suppose next that $M^\infty$ is noncompact. If $(M^\infty,g^\infty(t))$ is locally isometric to a round cylinder, then the point $x_j$ lies at the center of an $\varepsilon$-neck or an $\varepsilon$-quotient neck, contrary to our assumption. Thus, $(M^\infty,g^\infty(t))$ cannot be locally isometric to a round cylinder. Applying Theorem \ref{canonical.neighborhood.preparation.1} to $(M^\infty,g^\infty(t))$ (and with $\varepsilon$ replaced by $\frac{\varepsilon}{2}$), we conclude that there exists a neighborhood $B^\infty \subset M^\infty$ of the point $x_\infty$ such that $B_{g^\infty(0)}(x_\infty,C_1^{-1}) \subset B^\infty \subset B_{g^\infty(0)}(x_\infty,C_1)$ and $C_2^{-1} \leq \text{\rm scal}(x,0) \leq C_2$ for all $x \in B^\infty$. Moreover, $B^\infty$ is either an $\frac{\varepsilon}{2}$-neck or an $\frac{\varepsilon}{2}$-cap. Therefore, if $j$ is sufficiently large, then there exists a neighborhood $B^{(j)} \subset M^{(j)}$ of the point $x_j$ such that $B_{g^{(j)}(0)}(x_j,(2C_1)^{-1}) \subset B^{(j)} \subset B_{g^{(j)}(0)}(x_j,2C_1)$ and $(2C_2)^{-1} \leq \text{\rm scal}(x,0) \leq 2C_2$ for all $x \in B^{(j)}$. Moreover, $B^{(j)}$ is either an $\varepsilon$-neck or an $\varepsilon$-cap. This again contradicts the definition of $x_j$. \\

\section{A Canonical Neighborhood Theorem in higher dimensions}

\label{canonical.neighborhood.theorem}

In this section, we consider solutions to the Ricci flow starting from compact initial metrics. We assume that the solutions satisfies the following curvature pinching condition:

\begin{definition}
Let $f: [0,\infty) \to [0,\infty)$ be a concave and increasing function such that $\lim_{s \to \infty} \frac{f(s)}{s} = \frac{1}{n-2}$. We say that a Riemannian manifold has $(f,\theta)$-pinched curvature if the curvature tensor lies in the set 
\begin{align*} 
&\{R = S + H \owedge \text{\rm id}: S \in PIC2, \, \text{\rm Ric}_0(S)=0, \\ 
&\hspace{32mm} f(\text{\rm tr}(H)) \, \text{\rm id} - H \geq 0, \\ 
&\hspace{32mm} \text{\rm tr}(H) - \theta \, \text{\rm scal}(S) \geq 0\} 
\end{align*} 
at each point.
\end{definition}

The following is the analogue of Perelman's Canonical Neighborhood Theorem in dimension $3$:

\begin{theorem}[cf. G.~Perelman \cite{Perelman1}, Theorem 12.1]
Given a function $f$ as above and positive numbers $\theta$, $\kappa$, and $\varepsilon$, we can find a positive number $r_0$ such that the following holds: Let $(M,g(t))$, $t \in [0,T)$, be a compact solution to the Ricci flow which has $(f,\theta)$-pinched curvature and is $\kappa$-noncollapsed on scales less than $1$. Moreover, suppose that $M$ does not contain any non-trivial incompressible space forms $S^{n-1}/\Gamma$. Then for any point $(x_0,t_0)$ with $t_0 \geq 1$ and $Q := \text{\rm scal}(x_0,t_0) \geq r_0^{-2}$, the solution in $\{(x,t): d_{g(t_0)}(x_0,x) < \varepsilon^{-\frac{1}{2}} \, Q^{-\frac{1}{2}}, \, 0 \leq t_0-t \leq \varepsilon^{-1} \, Q^{-1}\}$ is, after scaling by the factor $Q$, $\varepsilon$-close to the corresponding subset of an ancient $\kappa_0$-solution with $\theta$-pinched curvature.
\end{theorem}

\textbf{Proof.} 
The proof is an adaptation of the argument in Section 12.1 of \cite{Perelman1}. We will follow the exposition in \cite{Chen-Zhu} and \cite{Kleiner-Lott}. Let $C_1=C_1(n,\theta,\varepsilon)$ denote the constant in Corollary \ref{canonical.neighborhood.preparation.2}. We define a constant $C_0=C_0(n,\theta,\varepsilon)$ by $C_0 := 4 \, \max \{C_1,\varepsilon^{-1}\}$.

Suppose that the assertion is false. We can find a sequence of Ricci flows $(M^{(j)},g^{(j)}(t))$ and a sequence of points $(x_j,t_j)$ in space-time with the following properties: 
\begin{itemize} 
\item[(i)] $(M^{(j)},g^{(j)}(t))$ does not contain any non-trivial incompressible $(n-1)$-dimensional space forms. 
\item[(ii)] $(M^{(j)},g^{(j)}(t))$ has $(f,\theta)$-pinched curvature.
\item[(iii)] $(M^{(j)},g^{(j)}(t))$ is $\kappa$-noncollapsed at all scales less than $1$. 
\item[(iv)] $t_j \geq \frac{1}{2}$ and $Q_j := \text{\rm scal}(x_j,t_j) \geq 2^j$. 
\item[(v)] After dilating by the factor $Q_j$, the solution in $\{(x,t): d_{g^{(j)}(t_j)}(x_j,x) < C_0^{\frac{1}{2}} \, Q_j^{-\frac{1}{2}}, \, 0 \leq t_j-t \leq C_0 \, Q_j^{-1}\}$ is not $\varepsilon$-close to the corresponding subset of any ancient $\kappa_0$-solution with $\theta$-pinched curvature.
\end{itemize}
By a standard point-picking argument, we can assume that $(x_j,t_j)$ in addition satisfies the following condition: 
\begin{itemize}
\item[(vi)] If $(\tilde{x},\tilde{t})$ is a point in space-time satisfying $\text{\rm scal}(\tilde{x},\tilde{t}) = : \tilde{Q} > 2Q_j$ and $0 \leq t_j-\tilde{t} \leq j \, Q_j^{-1}$, then the solution in $\{(x,t): d_{g^{(j)}(\tilde{t})}(\tilde{x},x) < C_0^{\frac{1}{2}} \, \tilde{Q}^{-\frac{1}{2}}, \, 0 \leq \tilde{t}-t \leq C_0 \, \tilde{Q}^{-1}\}$ is, after scaling by the factor $\tilde{Q}$, $\varepsilon$-close to the corresponding subset of an ancient $\kappa_0$-solution with $\theta$-pinched curvature.
\end{itemize}
Our goal is to show that, if we dilate the flow $(M^{(j)},g^{(j)}(t))$ around the point $(x_j,t_j)$ by the factor $Q_j$, the rescaled flows will converge in $C_{loc}^\infty$ to an ancient $\kappa_0$-solution with pinched curvature. To that end, we proceed in several steps: \\

\textit{Step 1:} The condition (vi) implies that $|D^m R(\tilde{x},\tilde{t})| \leq \eta \, \text{\rm scal}(\tilde{x},\tilde{t})^{\frac{m}{2}+1}$ whenever $m \in \{1,\hdots,8\}$, $\text{\rm scal}(\tilde{x},\tilde{t}) > 2Q_j$, and $0 \leq t_j-\tilde{t} \leq j \, Q_j^{-1}$. Here, $\eta$ is a large constant which is independent of $j$.

In particular, if $(\tilde{x},\tilde{t})$ is a point in space-time satisfying $0 \leq t_j - \tilde{t} \leq \frac{1}{2} \, j \, Q_j^{-1}$, then the gradient estimate implies that $\text{\rm scal}(x,t) \leq 4 \, (Q_j+\text{\rm scal}(\tilde{x},\tilde{t}))$ for all points $(x,t)$ satisfying $0 \leq \tilde{t}-t \leq c \, (Q_j+\text{\rm scal}(\tilde{x},\tilde{t}))^{-1}$ and $d_{g^{(j)}(\tilde{t})}(\tilde{x},x) \leq c \, (Q_j+\text{\rm scal}(\tilde{x},\tilde{t}))^{-\frac{1}{2}}$. Here, $c$ is a small positive constant which is independent of $j$. \\

\textit{Step 2:} For each $\rho \geq 0$, let $\mathbb{M}(\rho)$ be the smallest positive number such that 
\[Q_j^{-1} \, \text{\rm scal}(\tilde{x},t_j) \leq \mathbb{M}(\rho)\] 
for all integers $j \in \mathbb{N}$ and all points $\tilde{x} \in M^{(j)}$ satisfying $Q_j^{\frac{1}{2}} \, d_{g^{(j)}(t_j)}(x_j,\tilde{x}) \leq \rho$. If no such number exists, we put $\mathbb{M}(\rho) = \infty$. 

Using the local curvature bound in Step 1, we conclude that $\mathbb{M}(\rho) < \infty$ if $\rho>0$ is sufficiently small. Let 
\[\rho^* = \sup \{\rho \geq 0: \mathbb{M}(\rho) < \infty\}.\] 
We claim that $\rho^* = \infty$. We argue by contradiction, and assume that $\rho^* < \infty$. We can find a sequence of points $\tilde{x}_j \in M^{(j)}$ such that $\limsup_{j \to \infty} Q_j^{\frac{1}{2}} \, d_{g^{(j)}(t_j)}(x_j,\tilde{x}_j) \leq \rho^*$ and $\liminf_{j \to \infty} Q_j^{-1} \, \text{\rm scal}(\tilde{x}_j,t_j) = \infty$. Let $\gamma_j$ be a minimizing geodesic in $(M^{(j)},g^{(j)}(t_j))$ joining $x_j$ and $\tilde{x}_j$, and let $z_j$ be the point on $\gamma_j$ closest to $\tilde{x}_j$ with $\text{\rm scal}(z_j,t_j) = 4 \, Q_j$. We denote by $\beta_j$ the segment of $\gamma_j$ from $z_j$ to $\tilde{x}_j$.

We next dilate the ball $B_{g(t_j)}(x_j,\rho^* \, Q_j^{-\frac{1}{2}})$ by the factor $Q_j^{\frac{1}{2}}$. After passing to a subsequence, the dilated balls converge in $C_{loc}^\infty$ to an incomplete smooth manifold $(B^\infty,g^\infty)$. Moreover, the geodesics $\gamma_j$ and $\beta_j$ converge to  minimizing geodesics $\gamma_\infty$ and $\beta_\infty$ in $(B^\infty,g^\infty)$. Finally, the points $x_j$ and $z_j$ converge to points $x_\infty$ and $z_\infty$ in $M^\infty$. Using the gradient estimate $|D R| \leq \eta \, \text{\rm scal}^{\frac{3}{2}}$, we conclude that the curvature of $(B^\infty,g^\infty)$ must blow-up along $\beta_\infty$. Moreover, the curvature tensor of $(B^\infty,g^\infty)$ lies in the cone $\mathcal{C}_{1,\theta}$.

In view of statement (vi) above, each point on $\beta_j$ has a neighborhood of size $C_0 \, \text{\rm scal}^{-\frac{1}{2}}$ which is $\varepsilon$-close to an ancient $\kappa_0$-solution. Passing to the limit as $j \to \infty$, we conclude that each point $q \in \beta_\infty$ has a neighborhood of size $C_0 \, \text{\rm scal}_{g^\infty}(q)^{-\frac{1}{2}}$ which is $\varepsilon$-close to an ancient $\kappa_0$-solution. In particular, for each point $q \in \beta_\infty$, we have $C_0 \, \text{\rm scal}_{g^\infty}(q)^{-\frac{1}{2}} \leq \rho^* - d_{g^\infty}(x_\infty,q)$. Moreover, each point $q \in \beta_\infty$ has a neighborhood of size $\frac{1}{2} \, C_0 \, \text{\rm scal}_{g^\infty}(q)^{-\frac{1}{2}}$ which is $\varepsilon$-close to an ancient $\kappa_0$-solution. Using Corollary \ref{canonical.neighborhood.preparation.2} together with the fact that $C_0 \geq 4C_1$, we conclude that each point $q \in \beta_\infty$ has a canonical neighborhood $B$ which is either a $2\varepsilon$-neck; or a $2\varepsilon$-cap; or a closed manifold diffeomorphic to $S^n/\Gamma$; or a $2\varepsilon$-quotient neck. Let us consider the various possibilities: 
\begin{itemize}
\item If the canonical neighborhood of $q$ is a closed manifold, then the curvature of $(B^\infty,g^\infty)$ is bounded. This contradicts the fact that the curvature of $(B^\infty,g^\infty)$ blows up along $\beta_\infty$. Therefore, this case cannot occur.
\item If the canonical neighborhood of $q$ is a quotient neck, then $M^{(j)}$ contains a quotient neck if $j$ is sufficiently large. Theorem \ref{incompressibility.compact.case} then implies that $M^{(j)}$ contains a non-trivial incompressible $(n-1)$-dimensional space form for $j$ sufficiently large, contrary to our assumption. Hence, this case cannot occur. 
\item Finally, if $\rho^* - d_{g^\infty}(x_\infty,q)$ is sufficiently small and the canonical neighborhood of $q$ is a $2\varepsilon$-cap, then $\beta_\infty$ must enter and exit this cap, but this is impossible since $\beta_\infty$ is a minimizing geodesic. Consequently, this case cannot occur if $\rho^* - d_{g^\infty}(x_\infty,q)$ is suffificiently small.
\end{itemize}
To summarize, if $q \in \beta_\infty$ and $\rho^* - d_{g^\infty}(x_\infty,q)$ is sufficiently small, then $q$ has a canonical neighborhood which is a $2\varepsilon$-neck. Let $U$ denote the union of the canonical neighborhoods of all points $q \in \beta_\infty$, where $\rho^* - d_{g^\infty}(x_\infty,q)$ is sufficiently small. 

By work of Hamilton \cite{Hamilton5}, $U$ admits a foliation by a one-parameter family of constant mean curvature spheres $\Sigma_s$. We parametrize the surfaces $\Sigma_s$ so that the surfaces $\Sigma_s$ are defined for $s>0$ small enough, and move outward (away from $x_\infty$ and towards the end of the horn) as $s$ decreases towards $0$. Moreover, we can arrange that the lapse function $v: \Sigma_s \to \mathbb{R}$ has mean value $1$ for each $s>0$. Note that $\frac{1}{2} \leq v \leq 2$ on each leaf $\Sigma_s$. In particular, $\rho^* - d_{g^\infty}(x_\infty,q)$ is comparable to $s$ for each point $q \in \Sigma_s$.

Let $H(s)$ denote the mean curvature of $\Sigma_s$. Then 
\[-H'(s) = \Delta_{\Sigma_s} v + |A|^2 \, v + \text{\rm Ric}_{g^\infty}(\nu,\nu) \, v \geq \Delta_{\Sigma_s} v + \frac{1}{n-1} \, H(s)^2 \, v.\] 
We now take the mean value over $\Sigma_s$. Using the fact that $v$ has mean value $1$ and $\Delta_{\Sigma_s} v$ has mean $0$, we obtain 
\[-H'(s) \geq \frac{1}{n-1} \, H(s)^2.\] 
From this, we deduce that 
\[H(s) \leq \frac{n-1}{s}.\] 
Using again the fact that $v$ has mean $1$, we obtain 
\[\frac{d}{ds} \text{\rm area}_{g^\infty}(\Sigma_s) = H(s) \int_{\Sigma_s} v = H(s) \, \text{\rm area}_{g^\infty}(\Sigma_s) \leq \frac{n-1}{s} \, \text{\rm area}_{g^\infty}(\Sigma_s).\] 
Consequently, the function $s^{1-n} \, \text{\rm area}_{g^\infty}(\Sigma_s)$ is monotone decreasing in $s$. Moreover, if the function $s^{1-n} \, \text{\rm area}_{g^\infty}(\Sigma_s)$ is constant in $s$, then $H(s) = \frac{n-1}{s}$, $-H'(s) = \frac{1}{n-1} \, H(s)^2 + o(s^{-2})$, $|A|^2 = \frac{1}{n-1} \, H(s)^2$, and $\text{\rm Ric}_{g^\infty}(\nu,\nu)$; this, in turn, implies that $v$ is constant equal to $1$, and the manifold $U$ is a cone. (Note that the opening angle of the cone must be very small, as every point lies on a $2\varepsilon$-neck.) 

Since the function $s^{1-n} \, \text{\rm area}_{g^\infty}(\Sigma_s)$ is monotone decreasing, we obtain 
\[\liminf_{s \to 0} s^{1-n} \, \text{\rm area}_{g^\infty}(\Sigma_s) > 0\] 
or, equivalently, 
\[\limsup_{s \to 0} \sup_{q \in \Sigma_s} s^{-2} \, \text{\rm scal}_{g^\infty}(q) < \infty.\] 
On the other hand, since $C_0 \, \text{\rm scal}_{g^\infty}(q)^{-\frac{1}{2}} \leq \rho^* - d_{g^\infty}(x_\infty,q)$ for each point $q \in \beta_\infty$, we know that 
\[\liminf_{s \to 0} \sup_{q \in \Sigma_s} s^{-2} \, \text{\rm scal}_{g^\infty}(q) > 0\] 
or, equivalently, 
\[\limsup_{s \to 0} s^{1-n} \, \text{\rm area}_{g^\infty}(\Sigma_s) < \infty.\]
Therefore, the function $s^{1-n} \, \text{\rm area}_{g^\infty}(\Sigma_s)$ converges to a finite non-zero limit as $s \to \infty$. This gives $H(s)=\frac{n-1}{s} + o(s^{-1})$, $-H'(s) = \frac{1}{n-1} \, H(s)^2 + o(s^{-2})$, $|A|^2 = \frac{1}{n-1} \, H(s)^2$, and $\text{\rm Ric}_{g^\infty}(\nu,\nu) = o(s^{-2})$ on $\Sigma_s$. From this, we deduce that $v=1+o(1)$. Hence, if we dilate the manifold $U$ by the factor $s^{-1}$ around an arbitrary point on $\Sigma_s$ and send $s \to 0$, then the rescaled metrics will converge to a smooth limit metric which is a piece of a cone. (Again, the opening angle of the cone must be very small, as every point lies on a $2\varepsilon$-neck.) Using the local curvature estimate in Step 1, we can locally extend the metric backwards in time to a solution of the Ricci flow. To summarize, we obtain a (locally defined) solution to the Ricci flow whose curvature tensor lies in $\mathcal{C}_{1,\theta}$ and which, at the final time, is a piece of a cone. This contradicts Proposition \ref{splitting.1}. Thus, $\rho^*=\infty$. \\

\textit{Step 3:} We now dilate the manifold $(M^{(j)},g^{(j)}(t_j))$ around the point $x_j$ by the factor $Q_j$. Using the curvature bounds established in Steps 1 and 2 together with the $\kappa$-noncollapsing condition, we conclude that, after passing to a subsequence, the rescaled manifolds converge in the Cheeger-Gromov sense to a smooth limit manifold, which we denote by $(M^\infty,g^\infty)$. The condition (vi) implies that $|D^m R(\tilde{x})| \leq \eta \, \text{\rm scal}(\tilde{x})^{\frac{m}{2}+1}$ for all $m \in \{1,\hdots,8\}$ and all points $\tilde{x} \in M^\infty$ satisfying $\text{\rm scal}(\tilde{x}) \geq 4$. Moreover, since $(M^{(j)},g^{(j)}(t_j))$ has $(f,\theta)$-pinched curvature, the curvature tensor of $(M^\infty,g^\infty)$ lies in the cone $\mathcal{C}_{1,\theta}$.

We claim that the limit $(M^\infty,g^\infty)$ has bounded curvature. Indeed, if $(M^\infty,g^\infty)$ has unbounded curvature, then property (vi) above implies that $(M^\infty,g^\infty)$ contains a sequence of necks with radii converging to $0$, but this is impossible in a manifold with nonnegative sectional curvature.

\textit{Step 4:} We next show that $(M^\infty,g^\infty)$ can be extended backwards in time to an ancient solution. For each $\tau \geq 0$, let $\mathbb{L}(\tau)$ be the smallest positive number such that 
\[\limsup_{j \to \infty} Q_j^{-1} \, \text{\rm scal}(\tilde{x}_j,\tilde{t}_j) \leq \mathbb{L}(\tau)\] 
for every sequence of points $(\tilde{x}_j,\tilde{t}_j)$ in space-time satisfying $0 \leq t_j-\tilde{t}_j \leq \tau \, Q_j^{-1}$ and $\limsup_{j \to \infty} Q_j^{\frac{1}{2}} \, d_{g^{(j)}(\tilde{t}_j)}(x_j,\tilde{x}_j) < \infty$. If no such number exists, we put $\mathbb{L}(\tau) = \infty$. 

We have shown in Step 3 that the curvature of $(M^\infty,g^\infty)$ is bounded from above by some constant $\Lambda \geq 1$. This implies $\mathbb{L}(0) \leq \Lambda$. Using the local curvature bound in Step 1, we conclude that $\mathbb{L}(\tau) \leq 8(\Lambda+1)$ if $\tau>0$ is sufficiently small. Let 
\[\tau^* = \sup \{\tau \geq 0: \mathbb{L}(\tau) < \infty\}.\] 
We first show that $\sup_{\tau \in [0,\tau^*)} \mathbb{L}(\tau) \leq \Lambda$. Indeed, if $\sup_{\tau \in [0,\tau^*)} \mathbb{L}(\tau) > \Lambda$, then we dilate the flow $(M^{(j)},g^{(j)}(t))$ around the point $(x_j,t_j)$ by the factor $Q_j$. Using the fact that $\mathbb{L}(\tau)<\infty$ for each $\tau \in [0,\tau^*)$ together with the $\kappa$-noncollapsing condition, we conclude that the rescaled flows converge in the Cheeger-Gromov sense to a solution $(M^\infty,g^\infty(t))$ which is defined for $t \in (-\tau^*,0]$ and has bounded curvature for each $t$. Moreover, since $(M^{(j)},g^{(j)}(t))$ has $(f,\theta)$-pinched curvature, the curvature tensor of the limit flow $(M^\infty,g^\infty(t))$ lies in the cone $\mathcal{C}_{1,\theta}$. Using the Harnack inequality, we conclude that the curvature of the limit flow $(M^\infty,g^\infty(t))$ is bounded from above by $\Lambda$. This implies $\sup_{\tau \in [0,\tau^*)} \mathbb{L}(\tau) \leq \Lambda$. 

We claim that $\tau^* = \infty$. To prove this, we consider an arbitrary time $\hat{\tau} \in (0,\tau^*)$, and put $\hat{t}_j := t_j - \hat{\tau} \, Q_j^{-1}$. Using the local curvature estimate in Step 1 and the fact that $\mathbb{L}(\hat{\tau}) \leq \Lambda$, we conclude that there exists a positive number $c$ (independent of $j$ and $\hat{\tau}$) with the property that 
\[\limsup_{j \to \infty} Q_j^{-1} \, \text{\rm scal}(\tilde{x}_j,\tilde{t}_j) \leq 8(\Lambda+1)\] 
for every sequence of points $(\tilde{x}_j,\tilde{t}_j)$ in space-time satisfying $0 \leq \hat{t}_j-\tilde{t}_j \leq c \, Q_j^{-1}$ and $\limsup_{j \to \infty} Q_j^{\frac{1}{2}} \, d_{g^{(j)}(\hat{t}_j)}(x_j,\tilde{x}_j) < \infty$. This implies 
\[\limsup_{j \to \infty} Q_j^{-1} \, \text{\rm scal}(\tilde{x}_j,\tilde{t}_j) \leq 8(\Lambda+1)\] 
for every sequence of points $(\tilde{x}_j,\tilde{t}_j)$ in space-time satisfying $0 \leq \hat{t}_j-\tilde{t}_j \leq c \, Q_j^{-1}$ and $\limsup_{j \to \infty} Q_j^{\frac{1}{2}} \, d_{g^{(j)}(\tilde{t}_j)}(x_j,\tilde{x}_j) < \infty$. Thus, $\mathbb{L}(\hat{\tau}+c) \leq 8(\Lambda+1)$. In particular, $\hat{\tau}+c \leq \tau^*$ by definition of $\tau^*$. Since $\hat{\tau} \in (0,\tau^*)$ is arbitrary and $c>0$ is independent of $\hat{\tau}$, we conclude that $\tau^*=\infty$. \\

\textit{Step 5:} Finally, we dilate the flow $(M^{(j)},g^{(j)}(t))$ around the point $(x_j,t_j)$ by the factor $Q_j$. Using the fact that $\sup_{\tau \in [0,\infty)} \mathbb{L}(\tau) \leq \Lambda$ together with the $\kappa$-noncollapsing condition, we conclude that the rescaled flows converge in the Cheeger-Gromov sense to an ancient solution which has bounded curvature and is $\kappa$-noncollapsed. As above, the curvature tensor of the limit flow lies in the cone $\mathcal{C}_{1,\theta}$. Thus, the limit flow is an ancient $\kappa$-solution with $\theta$-pinched curvature. By Theorem \ref{universal.noncollapsing}, the limit flow is $\kappa_0$-noncollapsed. This contradicts statement (v) above. This completes the proof. \\

\begin{corollary}[cf. G.~Perelman \cite{Perelman1}, Theorem 12.1]
\label{high.curvature.regions}
Given a function $f$ as above and positive numbers $\theta$, $\kappa$, and $\varepsilon$, we can find a positive number $r_0$ such that the following holds: Let $(M,g(t))$, $t \in [0,T)$, be a compact solution to the Ricci flow which has $(f,\theta)$-pinched curvature and is $\kappa$-noncollapsed on scales less than $1$. Moreover, suppose that $M$ does not contain any non-trivial incompressible space forms $S^{n-1}/\Gamma$. Then for any point $(x_0,t_0)$ with $t_0 \geq 1$ and $Q := \text{\rm scal}(x_0,t_0) \geq r_0^{-2}$, there exists a neighborhood $B$ of $x_0$ such that $B_{g(t_0)}(x_0,(2C_1)^{-1} \, \text{\rm scal}(x_0,t_0)^{-\frac{1}{2}}) \subset B \subset B_{g(t_0)}(x_0,2C_1 \, \text{\rm scal}(x_0,t_0)^{-\frac{1}{2}})$ and $(2C_2)^{-1} \, \text{\rm scal}(x_0,t_0) \leq \text{\rm scal}(x,t_0) \leq 2C_2 \, \text{\rm scal}(x_0,t_0)$ for all $x \in B$. Finally, one of the following statements holds:
\begin{itemize}
\item $B$ is a $2\varepsilon$-neck.
\item $B$ is a $2\varepsilon$-cap in the sense that $B$ is diffeomorphic to a ball and the boundary $\partial B$ is the cross-sectional sphere of a $2\varepsilon$-neck. 
\item $B$ is a closed manifold diffeomorphic to $S^n/\Gamma$.
\end{itemize}
Here, $C_1=C_1(n,\theta,\varepsilon)$ and $C_2=C_2(n,\theta,\varepsilon)$ are the constants appearing in Corollary \ref{canonical.neighborhood.preparation.2}.
\end{corollary}

\section{The behavior of the flow at the first singular time}

\label{first.singular.time}

As in \cite{Perelman2}, we are now able to give a precise description of the behavior of the flow at the first singular time. Let us fix a compact initial manifold $(M,g_0)$ with the property that the curvature tensor of $(M,g_0)$ lies in the interior of the cone $\mathcal{C}_{2,0}$, and let $(M,g(t))$, $t \in [0,T)$, denote the unique maximal solution to the Ricci flow with initial metric $g_0$. We can find $\sigma_0 \in (1,2)$ and $\theta \in (0,\bar{\theta})$ such that the curvature tensor of $(M,g_0)$ lies in the interior of the cone $\mathcal{C}_{\sigma_0,\theta}$. By Theorem \ref{pinching.estimate}, there exists a concave and increasing function $f: [0,\infty) \to [0,\infty)$ such that $\lim_{s \to \infty} \frac{f(s)}{s} = \frac{1}{n-2}$; the manifold $(M,g_0)$ has $(f,\theta)$-pinched curvature; and $(f,\theta)$-pinching is preserved by the Ricci flow. In particular, $(M,g(t))$ has $(f,\theta)$-pinched curvature for all $t \in [0,T)$. Moreover, by work of Perelman \cite{Perelman1}, there exists a constant $\kappa>0$ such that $(M,g(t))$ is $\kappa$-noncollapsed for all $t \in [0,T)$. This allows us to apply Corollary \ref{high.curvature.regions}. In particular, we have $|D R| \leq \eta \, \text{\rm scal}^{\frac{3}{2}}$ and $|D^2 R| \leq \eta \, \text{\rm scal}^2$ at all points where the scalar curvature is sufficiently large. This implies that the set 
\[\Omega := \{x \in M: \limsup_{t \to T} \text{\rm scal}(x,t) < \infty\}\] 
is open. There are two possibilities now: 

The first possibility is that $\Omega$ is empty and the curvature becomes unbounded at each point on $M$. In this case, Corollary \ref{high.curvature.regions} implies that either $M$ is diffeomorphic to $S^n/\Gamma$, or else, for $t$ sufficiently close to $T$, every point $(x,t)$ lies on a $2\varepsilon$-neck or a $2\varepsilon$-cap. Thus, $M$ is diffeomorphic to to a quotient of $S^n$, or a quotient of $S^{n-1} \times \mathbb{R}^n$.

The second possibility is that $\Omega$ is non-empty. Since the derivatives of the curvature tensor can be controlled by suitable powers of the scalar curvature, the metrics $g(t)$ converge to a smooth limit metric $g(T)$ on $\Omega$. As in Perelman's paper \cite{Perelman2}, let 
\[\Omega_\rho := \{x \in M: \limsup_{t \to T} \text{\rm scal}(x,t) \leq \rho^{-2}\} = \{x \in \Omega: \text{\rm scal}(x,T) \leq \rho^{-2}\},\] 
where $\rho$ is chosen to be smaller than the scale $r_0$ in Corollary \ref{high.curvature.regions}. Then $\Omega_\rho$ is compact. By Corollary \ref{high.curvature.regions}, each point in $\Omega \setminus \Omega_\rho$ is contained in a subset of $\Omega$ which is either a $2\varepsilon$-neck or a $2\varepsilon$-cap or a closed manifold diffeomorphic to $S^n/\Gamma$. We can follow each $2\varepsilon$-neck to either side, until we reach the boundary of $\Omega$ (i.e. where the metric $g(T)$ becomes singular) or the boundary of $\Omega_\rho$ (i.e. where Corollary \ref{high.curvature.regions} can no longer be applied) or a $2\varepsilon$-cap. Therefore, each point in $\Omega \setminus \Omega_\rho$ is contained in a subset of $\Omega$ which falls into one of the following categories (in the terminology of \cite{Perelman2}): 
\begin{itemize}
\item an $2\varepsilon$-tube with boundary components in $\Omega_\rho$ 
\item an $2\varepsilon$-cap with boundary in $\Omega_\rho$
\item an $2\varepsilon$-horn with boundary in $\Omega_\rho$ 
\item a double $2\varepsilon$-horn
\item a capped $2\varepsilon$-horn
\item a closed manifold diffeomorphic to $S^n/\Gamma$
\end{itemize}
As in Perelman's paper \cite{Perelman2}, we perform surgery on every $2\varepsilon$-horn with boundary in $\Omega_\rho$. Moreover, we remove all double $2\varepsilon$-horns, all capped $2\varepsilon$-horns, and all closed manifolds diffeomorphic to $S^n/\Gamma$. The $2\varepsilon$-tubes and $2\varepsilon$-caps with boundary in $\Omega_\rho$ are left unchanged.

\begin{proposition}
\label{reconstructing.topology}
The pre-surgery manifold is a connected sum of the post-surgery manifold together with finitely many additional pieces. Each additional piece is diffeomorphic to a quotient of $S^n$ or to a compact quotient of $S^{n-1} \times \mathbb{R}$.
\end{proposition}

We next verify that our curvature pinching conditions are preserved under surgery. This can be viewed as the higher dimensional analogue of a theorem of R.~Hamilton (cf. \cite{Hamilton5}, Theorem D3.1):

\begin{proposition}
\label{pinching.preserved.under.surgery}
If a neck has $(f,\theta)$-pinched curvature right before surgery, then the surgically modified metric again has $(f,\theta)$-pinched curvature.
\end{proposition}

\textbf{Proof.} 
By scaling, we may assume that the radius of the neck is close to $1$. Let $g$ be a smooth metric on $S^{n-1} \times [-10,10]$ which is close to the standard metric and which has $(f,\theta)$-pinched curvature. Let $z$ denote the height function on $S^{n-1} \times [-10,10]$. As in \cite{Hamilton5}, the surgically modified metric is given by $\tilde{g} = e^{-2\varphi} \, g$, where $\varphi = e^{-\frac{1}{z}}$ for $z>0$ and $\varphi=0$ for $z \leq 0$. 

Let $\{e_1,\hdots,e_n\}$ denote a local orthonormal frame with respect to the metric $g$. If we put $\tilde{e}_i = e^\varphi e_i$, then $\{\tilde{e}_1,\hdots,\tilde{e}_n\}$ is an orthonormal frame with respect to the metric $\tilde{g}$. We will express geometric quantities associated with the metric $g$ relative to the frame $\{e_1,\hdots,e_n\}$, while geometric quantities associated with $\tilde{g}$ will be expressed in terms of $\{\tilde{e}_1,\hdots,\tilde{e}_n\}$. For example, we denote by $R_{ijkl} = R_g(e_i,e_j,e_k,e_l)$ the components of the Riemann curvature tensor of $g$, and by $\tilde{R}_{ijkl} = R_{\tilde{g}}(\tilde{e}_i,\tilde{e}_j,\tilde{e}_k,\tilde{e}_l)$ the components of the Riemann curvature tensor of $\tilde{g}$. With this understood, the standard formula for the change of the curvature tensor under a conformal change of the metric takes the form 
\[\tilde{R} = e^{2\varphi} \, R + e^{2\varphi} \, \Big ( D^2 \varphi + d\varphi \otimes d\varphi - \frac{1}{2} \, |d\varphi|^2 \, \text{\rm id} \Big ) \owedge \text{\rm id}\] 
(cf. \cite{Besse}, Theorem 1.159). Since $g$ has $(f,\theta)$-pinched curvature, the curvature tensor of $g$ can be written as $R = S + H \owedge \text{\rm id}$, where $S \in PIC2$, $\text{\rm Ric}_0(S)=0$, $f(\text{\rm tr}(H)) \, \text{\rm id} - H \geq 0$, and $\text{\rm tr}(H) - \theta \, \text{\rm scal}(S) \geq 0$. Consequently, the curvature tensor of $\tilde{g}$ can be written as $\tilde{R} = \tilde{S} + \tilde{H} \owedge \text{\rm id}$, where 
\[\tilde{S} = e^{2\varphi} \, S\] 
and 
\[\tilde{H} = e^{2\varphi} \, H + e^{2\varphi} \, \Big ( D^2 \varphi + d\varphi \otimes d\varphi - \frac{1}{2} \, |d\varphi|^2 \, \text{\rm id} \Big ).\] 
Clearly, $\tilde{S} \in PIC2$ and $\text{\rm Ric}_0(\tilde{S})=0$. We claim that $f(\text{\rm tr}(\tilde{H})) \, \text{\rm id} - \tilde{H} \geq 0$ and $\text{\rm tr}(\tilde{H}) - \theta \, \text{\rm scal}(\tilde{S}) \geq 0$, at least if $z>0$ is sufficiently small. To prove this, we observe that $D^2 \varphi = z^{-4} \, e^{-\frac{1}{z}} \, dz \otimes dz + o(z^{-4} \, e^{-\frac{1}{z}})$. Hence, if $z>0$ is sufficiently small, we have 
\[\text{\rm tr}(\tilde{H}) - \text{\rm tr}(H) \geq (1-c) \, z^{-4} \, e^{-\frac{1}{z}},\] 
where $c$ is a positive constant that can be made arbitrarily small. Since the function $f(s)-\frac{s}{n-2}$ is monotone increasing in $s$, it follows that 
\[f(\text{\rm tr}(\tilde{H})) - f(\text{\rm tr}(H)) \geq \frac{1}{n-2} \, (\text{\rm tr}(\tilde{H}) - \text{\rm tr}(H)) \geq \frac{1-c}{n-2} \, z^{-4} \, e^{-\frac{1}{z}}.\] 
Since the eigenvector corresponding to the smallest eigenvalue of $H$ is nearly orthogonal to the cross-section of the neck, we obtain 
\[\lambda_n(\tilde{H}) - \lambda_n(H) \leq c \, z^{-4} \, e^{-\frac{1}{z}}\] 
if $z>0$ is sufficiently small. As above, $c$ is a positive constant that can be made arbitrarily small. Putting these facts together, we conclude that 
\[f(\text{\rm tr}(\tilde{H})) - \lambda_n(\tilde{H}) > f(\text{\rm tr}(H)) - \lambda_n(H) \geq 0\] 
if $z>0$ is sufficiently small. Moreover, the inequality $\text{\rm tr}(\tilde{H}) > e^{2\varphi} \, \text{\rm tr}(H)$ gives 
\[\text{\rm tr}(\tilde{H}) - \theta \, \text{\rm scal}(\tilde{S}) > e^{2\varphi} \, (\text{\rm tr}(H) - \theta \, \text{\rm scal}(S)) \geq 0\] 
if $z>0$ is sufficiently small. This shows that the surgically modified metric has $(f,\theta)$-pinched curvature when $z>0$ is sufficiently small. Finally, if $z$ is bounded away from $0$, then it is easy to see that the curvature tensor of the surgically modified metric lies in the cone $\mathcal{C}_{1,\theta}$. Therefore, the surgically modified metric has $(f,\theta)$-pinched curvature everywhere. \\

We now consider the surgically modified manifold, and evolve it again by the Ricci flow. The arguments in Section 5 of \cite{Chen-Zhu} (which are based on Section 5 in Perelman's paper \cite{Perelman2}) show that the surgery cutoff can be chosen in such a way that the $\kappa$-noncollapsing property and the Canonical Neighborhood Theorem (Corollary \ref{high.curvature.regions}) hold for the surgically modified flow. From this, Theorem \ref{surgically.modified.flow} follows.

\appendix 

\section{A higher-dimensional version of Theorem C4.1 in Hamilton's paper \cite{Hamilton5}}

In this final section, we collect some auxiliary results which are needed to rule out quotient necks. Recall that a quotient neck is modeled on a noncompact quotient of the cylinder $S^{n-1} \times \mathbb{R}$. Let $\Gamma \subset \text{\rm Isom}(S^{n-1} \times \mathbb{R})$ be a discrete group such that the quotient $(S^{n-1} \times \mathbb{R})/\Gamma$ is a noncompact smooth manifold, and let $\Gamma_1$ denote the image of $\Gamma$ under the canonical projection $\text{\rm Isom}(S^{n-1} \times \mathbb{R}) \to \text{\rm Isom}(\mathbb{R})$. Since the quotient $(S^{n-1} \times \mathbb{R})/\Gamma$ is noncompact, the group $\Gamma_1$ is either trivial, or it consists of the identity and a reflection. Without loss of generality, we may assume that $\Gamma_1$ fixes $0$, so that the slice $S^{n-1} \times \{0\}$ is invariant under $\Gamma$. 

\begin{theorem}[cf. R.~Hamilton \cite{Hamilton5}, Theorem C4.1] 
\label{incompressibility.compact.case}
Let $M$ be a compact manifold with positive isotropic curvature, and let $\Gamma$ be a discrete group as above. Suppose that $(S^{n-1} \times [-L,L])/\Gamma$ can be embedded into $M$ such that the induced metric is sufficiently close to the standard metric on $(S^{n-1} \times [-L,L])/\Gamma$ in the $C^3$-norm. Then the map $\pi_1((S^{n-1} \times \{0\})/\Gamma) \to \pi_1(M)$ is injective. 
\end{theorem} 

\textbf{Proof.} 
\textit{Case 1:} Assume first that $\Gamma_1$ is trivial. In this case, each slice $S^{n-1} \times \{z\}$ is invariant under the action of $\Gamma$. Suppose that the map $\pi_1((S^{n-1} \times \{0\})/\Gamma) \to \pi_1(M)$ is not injective. Then there exists a disk $D \subset M$ such that the boundary $\partial D$ is contained in $(S^{n-1} \times \{0\})/\Gamma$ and $\partial D$ represents a non-trivial element of $\pi_1((S^{n-1} \times \{0\})/\Gamma)$. Without loss of generality, we may assume that $D$ intersects $(S^{n-1} \times \{0\})/\Gamma$ transversally. The intersection $D \cap ((S^{n-1} \times \{0\})/\Gamma)$ can be written as a union of $\partial D$ with finitely many closed curves $\gamma_1,\hdots,\gamma_l$ ($l \geq 0$). Let us choose a disk $D$ for which $l$ is as small as possible.

We claim that $l = 0$. Suppose that $l \geq 1$. The curve $\gamma_1$ bounds a disk $\tilde{D} \subset D$. The intersection $\tilde{D} \cap ((S^{n-1} \times \{0\})/\Gamma)$ has fewer connected components than $D \cap ((S^{n-1} \times \{0\})/\Gamma)$. By minimality of $l$, the boundary $\partial \tilde{D}=\gamma_1$ must be homotopically trivial in $(S^{n-1} \times \{0\})/\Gamma$. Hence, there exists a closed curve $\hat{\gamma}_1 \subset D$ which encloses $\gamma_1$ and which bounds a disk in $M \setminus ((S^{n-1} \times \{0\})/\Gamma)$. This allows us to construct a disk $\hat{D} \subset M$ such that $\partial \hat{D} = \partial D$, $\hat{D}$ intersects $(S^{n-1} \times \{0\})/\Gamma$ transversally, and $\hat{D} \cap ((S^{n-1} \times \{0\})/\Gamma)$ has fewer connected components than $D \cap ((S^{n-1} \times \{0\})/\Gamma)$. This contradicts the minimality of $l$. Therefore, $l=0$. This means that $D$ does not intersect $(S^{n-1} \times \{0\})/\Gamma$ away from the boundary $\partial D$.

We now cut $M$ open along the center slice $(S^{n-1} \times \{0\})/\Gamma$ to get a (possibly disconnected) manifold $\hat{M}$ with boundary $\partial \hat{M}$ consisting of two copies of $(S^{n-1} \times \{0\})/\Gamma$. We can perturb the metric so that $\hat{M}$ has positive isotropic curvature and $\partial \hat{M}$ is convex. By work of Fraser (cf. \cite{Fraser}, Theorem 2.1b, Theorem 3.1, and Theorem 3.2), $\pi_2(\hat{M},\partial \hat{M}) = 0$. On the other hand, the disk $D$ constructed above represents a non-trivial element of $\pi_2(\hat{M},\partial \hat{M})$. This is a contradiction. 

\textit{Case 2:} Assume next that $n$ is even and $\Gamma_1$ has order $2$. Let $\tilde{\Gamma}$ denote the subgroup consisting of all elements of $\Gamma$ which map each slice $S^{n-1} \times \{z\}$ to itself. Since $n$ is even, every fixed point free isometry of $S^{n-1}$ preserves orientation. Hence, $\tilde{\Gamma}$ can be characterized as the subgroup consisting of all elements of $\Gamma$ which preserve orientation. 

Let $\tilde{M}$ denote the orientable double cover of $M$. Then $(S^{n-1} \times [-L,L])/\tilde{\Gamma}$ can be embedded into $\tilde{M}$ such that the induced metric is close to the standard metric on $(S^{n-1} \times [-L,L])/\tilde{\Gamma}$ in the $C^3$-norm. The arguments in Case 1 imply that the map $\pi_1((S^{n-1} \times \{0\})/\tilde{\Gamma}) \to \pi_1(\tilde{M})$ is injective. This easily implies that the map $\pi_1((S^{n-1} \times \{0\})/\Gamma) \to \pi_1(M)$ is injective. 


\textit{Case 3:} Assume finally that $n$ is odd and $\Gamma_1$ has order $2$. In this case, $\Gamma$ has order $2$ and consists of the identity and the reflection which sends $(p,z) \in S^{n-1} \times \mathbb{R}$ to $(-p,-z) \in S^{n-1} \times \mathbb{R}$. Let $\alpha$ be a closed curve in $(S^{n-1} \times \{0\})/\Gamma$ which represents the non-trivial element of $\pi_1((S^{n-1} \times \{0\})/\Gamma)$. We can deform $\alpha$ to a closed curve which intersects the center slice $(S^{n-1} \times \{0\})/\Gamma$ transversally at exactly one point. This shows that $\alpha$ represents a non-trivial element of $\pi_1(M)$. \\

\begin{theorem}[cf. R.~Hamilton \cite{Hamilton5}, Theorem C4.1] 
\label{incompressibility.noncompact.case}
Let $M$ be a complete noncompact manifold with curvature tensor in the interior of the PIC2 cone, and let $\Gamma$ be a discrete group as above. Suppose that $(S^{n-1} \times [-L,L])/\Gamma$ can be embedded into $M$ such that the induced metric is sufficiently close to the standard metric on $(S^{n-1} \times [-L,L])/\Gamma$ in the $C^3$-norm. Then the map $\pi_1((S^{n-1} \times \{0\})/\Gamma) \to \pi_1(M)$ is injective. 
\end{theorem} 

\textbf{Proof.} 
Again, we first consider the case when $\Gamma_1$ is trivial. If the map $\pi_1((S^{n-1} \times \{0\})/\Gamma) \to \pi_1(M)$ is not injective, we find a disk $D \subset M$ such that $\partial D$ is contained in $(S^{n-1} \times \{0\})/\Gamma$ and $\partial D$ represents a non-trivial element of $\pi_1((S^{n-1} \times \{0\})/\Gamma)$. As above, we can arrange that $D$ meets $(S^{n-1} \times \{0\})/\Gamma$ transversally along $\partial D$ and $D$ does not intersect $(S^{n-1} \times \{0\})/\Gamma$ away from the boundary $\partial D$. By \cite{Greene-Wu}, $M$ admits a smooth, strictly convex, exhaustion function. Hence, we can find a smooth, strictly convex domain $U \subset M$ which contains the neck $(S^{n-1} \times [-L,L])/\Gamma$ and the disk $D$. We then cut $U$ open along the center slice $(S^{n-1} \times \{0\})/\Gamma$. This gives a manifold $\hat{U}$ with boundary $\partial \bar{U}$ given by a union of $\partial U$ and two copies of $(S^{n-1} \times \{0\})/\Gamma$. We can perturb the metric so that $\hat{U}$ has positive isotropic curvature and $\partial \hat{U}$ is convex. Again, results in \cite{Fraser} imply $\pi_2(\hat{U},\partial \hat{U}) = 0$. This contradicts the fact that $D$ represents a non-trivial element of $\pi_2(\hat{U},\partial \hat{U})$.

The case when $\Gamma_1$ has order $2$ proceeds the same way as before.


\begin{thebibliography}{99}
\bibitem{Besse}
A.L.~Besse, \textit{Einstein manifolds,} Ergebnisse der Mathematik und ihrer Grenzgebiete, vol. 10, Springer-Verlag, Berlin (1987)

\bibitem{Bohm-Wilking} 
C.~B\"ohm and B.~Wilking, \textit{Manifolds with positive curvature operator are space forms,} Ann. of Math. 167, 1079--1097 (2008)

\bibitem{Bony} 
J.M.~Bony, \textit{Principe du maximum, in\'egalit\'e de Harnack et unicit\'e du probl\`eme de Cauchy pour les op\'erateurs elliptiques d\'eg\'en\'er\'es,} Ann. Inst. Fourier (Grenoble) 19, 277--304 (1969)

\bibitem{Brendle1}
S.~Brendle, \textit{A general convergence result for the Ricci flow in higher dimensions,} Duke Math. J. 145, 585--601 (2008)

\bibitem{Brendle2}
S.~Brendle, \textit{A generalization of Hamilton's differential Harnack inequality for the Ricci flow,} J. Diff. Geom. 82, 207--227 (2009)

\bibitem{Brendle-book}
S.~Brendle, \textit{Ricci Flow and the Sphere Theorem,} Graduate Studies in Mathematics, vol. 111, American Mathematical Society (2010)

\bibitem{Brendle-Schoen} 
S.~Brendle and R.~Schoen, \textit{Manifolds with $1/4$-pinched curvature are space forms,} J. Amer. Math. Soc. 22, 287--307 (2009)

\bibitem{Cheeger-Gromoll1} 
J.~Cheeger and D.~Gromoll, \textit{The splitting theorem for manifolds of nonnegative Ricci curvature,} J. Diff. Geom. 6, 119--128 (1971)

\bibitem{Cheeger-Gromoll2}
J.~Cheeger and D.~Gromoll, \textit{On the structure of complete manifolds of nonnegative curvature,} Ann. of Math. 96, 413--443 (1972)

\bibitem{Chen-Zhu} 
B.~Chen and X.~Zhu, \textit{Ricci flow with surgery on four-manifolds with positive isotropic curvature,} J. Diff. Geom. 74, 177--264 (2006) 

\bibitem{Chen-Tang-Zhu}
B.~Chen, S.~Tang, and X.~Zhu, \textit{Complete classification of compact four-manifolds with positive isotropic curvature,} J. Diff. Geom. 91, 41--80 (2012)

\bibitem{Fraser}
A.M.~Fraser, \textit{Minimal disks and two-convex hypersurfaces,} Amer. J. Math. 124, 483--493 (2002)

\bibitem{Greene-Wu}
R.E.~Greene and H.~Wu, \textit{$C^\infty$ convex functions and manifolds of positive curvature,} Acta Math. 137, 209--245 (1976)

\bibitem{Hamilton1}
R.~Hamilton, \textit{Three-manifolds with positive Ricci curvature,} J. Diff. Geom. 17, 255--306 (1982)

\bibitem{Hamilton2} 
R.~Hamilton, \textit{Four-manifolds with positive curvature operator,} J. Diff. Geom. 24, 153--179 (1986)

\bibitem{Hamilton3} 
R.~Hamilton, \textit{The Harnack estimate for the Ricci flow,} J. Diff. Geom. 37, 225--243 (1993)

\bibitem{Hamilton4}
R.~Hamilton, \textit{The formation of singularities in the Ricci flow,} Surveys in Differential Geometry, vol. II, 7--136, International Press, Somerville MA (1995)

\bibitem{Hamilton5}
R.~Hamilton, \textit{Four-manifolds with positive isotropic curvature,} Comm. Anal. Geom. 5, 1--92 (1997)

\bibitem{Huisken} 
G.~Huisken, \textit{Ricci deformation of the metric on a Riemannian manifold,} J. Diff. Geom. 21, 47--62 (1985)

\bibitem{Ivey}
T.~Ivey, \textit{New examples of complete Ricci solitons,} Proc. Amer. Math. Soc. 122, 241--245 (1994)

\bibitem{Kleiner-Lott}
B.~Kleiner and J.~Lott, \textit{Notes on Perelman's papers,} Geom. Topol. 12, 2587--2855 (2008)

\bibitem{Margerin1}
C.~Margerin, \textit{Pointwise pinched manifolds are space forms,} Geometric measure theory and the calculus of variations (Arcata 1984), 307--328, Proc. Sympos. Pure Math. 44, Amer. Math. Soc., Providence RI (1986)

\bibitem{Margerin2}
C.~Margerin, \textit{A sharp characterization of the smooth $4$-sphere in curvature terms,} Comm. Anal. Geom. 6, 21--65 (1998)

\bibitem{Margerin3}
C.~Margerin, \textit{D\'eformations de structures Riemanniennes,} unpublished manuscript

\bibitem{Micallef-Moore}
M.~Micallef and J.D.~Moore, \textit{Minimal two-spheres and the topology of manifolds with positive curvature on totally isotropic two-planes,} Ann. of Math. 127, 199--227 (1988)

\bibitem{Nishikawa}
S.~Nishikawa, \textit{Deformation of Riemannian metrics and manifolds with bounded curvature ratios,} Geometric measure theory and the calculus of variations (Arcata 1984), 343--352, Proc. Sympos. Pure Math. 44, Amer. Math. Soc., Providence RI (1986)

\bibitem{Perelman1} 
G.~Perelman, \textit{The entropy formula for the Ricci flow and its geometric applications,} arxiv:0211159

\bibitem{Perelman2}
G.~Perelman, \textit{Ricci flow with surgery on three-manifolds,} arxiv:0303109

\bibitem{Perelman3}
G.~Perelman, \textit{Finite extinction time for solutions to the Ricci flow on certain three-manifolds,} arxiv:0307245
\end{thebibliography}
\end{document}